\documentclass[11pt]{amsart}
\usepackage{amsmath}
\usepackage{amsxtra}
\usepackage{amscd}
\usepackage{amsthm}
\usepackage{amsfonts}
\usepackage{amssymb}
\usepackage{eucal}

\newtheorem{thm}{Theorem}
\newtheorem{lem}{Lemma}
\newtheorem{cor}{Corollary}
\newtheorem{prop}{Proposition}
\newtheorem{conj}{Conjecture}

\theoremstyle{definition}

\theoremstyle{remark}
\newtheorem{rem}{Remark}[section]

\numberwithin{equation}{section}

\begin{document}

\newcommand{\thmref}[1]{Theorem~\ref{#1}}
\newcommand{\conjref}[1]{Conjecture~\ref{#1}}
\newcommand{\secref}[1]{Sect.~\ref{#1}}
\newcommand{\lemref}[1]{Lemma~\ref{#1}}
\newcommand{\propref}[1]{Proposition~\ref{#1}}
\newcommand{\corref}[1]{Corollary~\ref{#1}}
\newcommand{\remref}[1]{Remark~\ref{#1}}
\newcommand{\nc}{\newcommand}
\nc{\on}{\operatorname}
\nc{\ch}{\mbox{ch}}
\nc{\Z}{{\mathbb Z}}
\nc{\C}{{\mathbb C}}
\nc{\cond}{|\,}
\nc{\bib}{\bibitem}
\nc{\pone}{\Pro^1}
\nc{\pa}{\partial}
\nc{\F}{{\mc F}}
\nc{\arr}{\rightarrow}
\nc{\larr}{\longrightarrow}
\nc{\al}{\alpha}
\nc{\ri}{\rangle}
\nc{\lef}{\langle}
\nc{\W}{{\mc W}}
\nc{\gam}{\ol{\gamma}}
\nc{\Q}{\ol{Q}}
\nc{\q}{\widetilde{Q}}
\nc{\la}{\lambda}
\nc{\ep}{\epsilon}
\nc{\su}{\widehat{\sw}_2}
\nc{\g}{\ol{{\mf g}}}
\nc{\hh}{{\mf h}}
\nc{\h}{\ol{\mf h}}
\nc{\n}{\ol{\mf n}}
\nc{\ab}{{\mf a}}
\nc{\f}{\widehat{{\mc F}}}
\nc{\laa}{\ol{\lambda}}
\nc{\De}{\Delta}
\nc{\G}{{\mf g}}
\nc{\Li}{{\mc L}}
\nc{\La}{\Lambda}
\nc{\is}{{\mathbf i}}
\nc{\V}{\widetilde{V}}
\nc{\M}{\widetilde{M}}
\nc{\js}{{\mathbf j}}
\nc{\sw}{{\mf s}{\mf l}}
\nc{\one}{{\mathbf 1}}
\nc{\N}{{\mf n}}
\nc{\HH}{{\mc H}}
\nc{\GG}{\widehat{G}}
\nc{\w}{\overline{W}}
\nc{\lo}{_{\on{loc}}}
\nc{\cc}{\widehat{C}}
\nc{\wt}{\widetilde}
\nc{\wh}{\widehat}
\nc{\bi}{\bibitem}

\nc{\mc}{\mathcal}
\nc{\mf}{\mathfrak}
\nc{\ol}{\overline}
\nc{\el}{\ell}
\nc{\WW}{\mathbf W}
\nc{\NN}{\ol{N}}
\nc{\paa}{\wt{\pa}}
\nc{\pq}{\wt{p} \; \wt{q}}
\nc{\ma}{\on{max}}
\nc{\rh}{\rho^\vee}
\nc{\rhh}{\rh_{-1}}

\title{Integrable hierarchies and Wakimoto modules}

\author[Boris Feigin]{Boris Feigin$^1$}\thanks{$^1$Partially supported
by the grants RFBR 99-01-01169, INTAS-OPEN-97-1312} \address{Landau
Institute for Theoretical Physics, Kosygina St 2, Moscow 117940,
Russia}

\author[Edward Frenkel]{Edward Frenkel$^2$}\thanks{$^2$Partially
supported by the Packard Foundation and the NSF} \address{Department
of Mathematics, University of California, Berkeley, CA 94720, USA}

\date{February 1999}

\dedicatory{To our teacher D.B. Fuchs on his sixtieth birthday}

\maketitle

\section{Introduction}

In our papers \cite{ff:im,ff:kdv} we proposed a new approach to
integrable hierarchies of soliton equations and their quantum
deformations. We have applied this approach to the Toda field theories
and the generalized KdV and modified KdV (mKdV) hierarchies. In this
paper we apply our approach to the Ablowitz-Kaup-Newell-Segur (AKNS)
hierarchy \cite{AKNS} and its generalizations. In particular, we show
that the free field (Wakimoto) realization of an affine algebra
\cite{Wak,ff:usp} naturally appears in the context of the generalized
AKNS hierarchies. This is analogous to the appearance of the free field
(quantum Miura) realization of a $\W$--algebra in the context of the
generalized KdV equations. As an application, we give here a new proof
of the existence of the Wakimoto realization.

We also conjecture that all integrals of motion of the generalized
AKNS equation can be quantized. In the case of $\widehat{\sw}_2$ the
corresponding quantum integrals of motion can be viewed as integrals
of motion of a thermal perturbation of the parafermionic conformal
field theory \cite{fateev}. Thus we expect that this deformation, and
analogous deformations for arbitrary affine algebras, are integrable,
in the sense of Zamolodchikov \cite{zam}.

Let us first recall the main steps of our analysis of the Toda field
theories from \cite{ff:im,ff:kdv}.

\subsection{Overview of the previous work} Let $\G$ be an affine
algebra, twisted or untwisted, and $\g$ be the finite-dimensional
simple Lie algebra, whose Dynkin diagram is obtained by deleting the
$0$th node of the Dynkin diagram of $\G$ \cite{K}. Let $\h$ be the
Cartan subalgebra of $\g$ and $\hh$ be the corresponding Heisenberg
algebra. For $x \in \hh, n \in \Z$, denote $x(n) = x \otimes
t^n$. Then $\h$ has generators $b_i(n) = \al_i(n), i=1,\ldots,\el,
n\in\Z$ (where $\al_i$'s are the simple roots of $\g$), and $\one$,
with the commutation relations
\begin{equation}    \label{heisrel}
[b_i(n),b_j(m)] = n (\al_i,\al_j) \delta_{n,-m} \one, \quad \quad
[\one,b_i(n)] = 0.
\end{equation}

There is a family of Fock representations $\pi^\nu_\la, \nu \in \C, \la \in
\h^*$, of $\hh$, which are generated by vectors $v_\la$, such that
$$b_i(n) v_\la = 0, n>0, \quad
\quad b_i(0) v_\la = (\la,\al_i) v_\la, \quad \quad \one v_\la = \nu
v_\la.$$ Introduce a derivation $\pa$ on $\pi^\nu_\la$, such that
$[\pa,b_i(n)] = -n b_i(n-1)$, $\pa \cdot v_\la =
\la(-1) v_\la$, where $\la \in \h \simeq \h^*$, and so $\la(-1) \in \hh$.

The module $\pi^\nu_0$ carries a vertex operator algebra (VOA)
structure described in \cite{ff:im} (for the definition of VOA, see
\cite{B,FLM}, or Sect. 4 of \cite{ff:im}).

Consider the space $\pi_\la[z,z^{-1}]$, and extend the action of
$\pa$ to it by the formula $\pa \otimes 1 + 1 \otimes
\pa_z$. Denote by $\F^\nu_\la$ the quotient of $\pi^\nu_\la \otimes
\C[z,z^{-1}]$ by the image of $\pa$ (total derivatives) and constants, if
$\la = 0$. The space $\F^\nu_0$ is a Lie algebra, the {\em local
completion} of the universal enveloping algebra of $\hh$
\cite{ff:gd}. Elements of $\F^\nu_0$ can be interpreted as Fourier
components of currents of the VOA $\pi^\nu_0$ \cite{ff:im}.

For any $\gamma \in \h^*$ we can define the bosonic vertex operator
\begin{equation}    \label{vertex}
V^\nu_\gamma(z) = \sum_{n\in\Z} V^\nu_\gamma(n) z^{\nu(\gamma,\la)-n} =
\end{equation}
$$T_\gamma z^{\nu(\gamma,\la)} \exp \left( -\sum_{n<0} \frac{\gamma(n)
z^{-n}}{n} \right) \exp \left( -\sum_{n>0} \frac{\gamma(n) z^{-n}}{n}
\right),$$
where $\gamma \in {\mf h}^* \simeq {\mf h}$ and $T_\gamma:
\pi^\nu_\la \arr \pi^\nu_{\la+\gamma}$ is such that $T_\gamma \cdot v_\la =
v_{\la+\gamma}$ and $[T^\nu_\gamma,b_i(n)]=0, n<0$. Thus, $V^\nu_\gamma(n),
n \in \Z$, are well-defined linear operators $\pi^\nu_\la \arr
\pi^\nu_{\la+\gamma}$.

Now introduce operators $$\q^\nu_i = V^\nu_{-\al_i}(1) = \int
V^\nu_{-\al_i}(z) dz: \pi_0 \arr \pi_{-\al_i}, \quad \quad
i=1,\ldots,\el,$$ where $\al_i, i=1,\ldots,\el$, are the simple roots
of $\g$. These operators are called the {\em screening
operators}. They commute with the action of $\pa$ and hence give rise
to operators $\Q^\nu_i: \F^\nu_0 \arr \F^\nu_{-\al_i},
i=1,\ldots,\el$.

The VOA $\WW_\nu(\g)$ is defined as a vertex operator subalgebra of
$\pi^\nu_0$: $$\WW_\nu(\g) = \bigcap_{i=1}^\el \on{Ker}_{\pi^\nu_0}
\q^\nu_i.$$

In \cite{ff:im} we proved that for generic $\nu$, $\WW_\nu(\g)$ is
finitely and freely generated in the following sense. There exist
elements $W^\nu_i$ in $\WW_\nu(\g)$ of degrees $d_i+1, i=1,\ldots,\el$,
where $d_i$ is the $i$th exponent of $\g$, such that $\WW_\nu(\g)$ has
a linear basis of lexicographically ordered monomials in the Fourier
components $W^\nu_i(n_i), 1\leq i\leq l,n_i<-d_i,$ of the currents
$Y(W^\nu_i,z) = \sum_{n\in\Z} W^\nu_i(n) z^{-n-d_i-1}$.

The Lie algebra $I_\nu(\g)$ of quantum integrals of motion of the Toda
theory of $\g$ (the $\W$--algebra of $\g$) is defined as $$\W_\nu(\g) =
\bigcap_{i=1}^\el \on{Ker}_{\F^\nu_0} \Q^\nu_i.$$

We proved in \cite{ff:im} that $I_\nu(\g)$ is isomorphic to the
quotient of $\WW_\nu(\g)[z,z^{-1}]$ by the total
derivatives and constants. In other words, it consists of all Fourier
components of currents defined by the VOA $\WW_\nu(\g)$. The Lie
algebra $I_\nu(\g)$ is called the {\em quantum $\W$--algebra}
associated to $\g$ and is denoted by $\W_\nu(\g)$.

Now let $\al_0 = - 1/a_0 \sum_{i=1}^\el a_i \al_i$ be the element of
$\h^*$, corresponding to the affine root of $\G$; here $a_i,
i=1,\ldots,\el$, are the labels of the Dynkin diagram of $\G$
\cite{K}. We can define the screening operators corresponding to the
$0$th root of the affine algebra $\G$, $\q^\nu_0: \pi^\nu_0 \arr
\pi^\nu_{-\al_0}$ and $\Q^\nu_0: \F^\nu_0 \arr \F^\nu_{-\al_0}$ in the
same way as the operators $\Q^\nu_i$ and $\q^\nu_i, i=1,\ldots,\el$,
which correspond to the simple roots of $\g$. The space $I_\nu(\G)$ of
{\em quantum integrals of motion} of the affine Toda theory associated
to $\G$ is $$I_\nu(\G) = \bigcap_{i=0}^\el \on{Ker}_{\F^\nu_0}
\Q^\nu_i.$$ Clearly, $I_\nu(\G)$ is a Lie subalgebra of
$\W_\nu(\g)$. Elements of the space $I_\nu(\G)$ can be viewed as
integrals of motion of a certain perturbation of the conformal field
theory with the $\WW_\nu(\g)$--symmetry (see
\cite{zam,EY,HM,ff:toda}). In \cite{ff:im} we proved that for generic
$\nu$ this is a commutative Lie algebra linearly spanned by elements
of degrees equal to the exponents of $\G$ modulo the Coxeter
number.\footnote{The intersection of kernels of the operators
$\q^\nu_i, i=0,\ldots,\el$, is one-dimensional and therefore not
interesting}

In order to prove the above mentioned results we interpreted the spaces
$\W_\nu(\g)$ and $I_\nu(\G)$ as cohomologies of certain complexes. These
complexes were constructed using the Bernstein-Gelfand-Gelfand (BGG)
resolution of the trivial representation of the quantum group $U_q(\g)$ and
$U_q(\G)$, respectively, where $q=\exp (\pi i \nu)$. The construction was
based on the fact that in a certain sense the operators $\Q^\nu_i$ satisfy
the defining relations of the nilpotent subalgebra of a quantum group --
the $q$--Serre relations \cite{bmp}.

We were able to compute the cohomologies of our complexes for generic
$\nu$ by computing them in the {\em classical limit} $\nu \arr 0$
\cite{ff:im}. In this limit, the spaces $I_\nu(\G)$ and $I_\nu(\g)$
become the spaces $I_0(\G)$ and $I_0(\g)$ of local integrals of motion
of the Toda equation associated to $\G$ and $\g$, respectively.

In the classical limit $\nu \arr 0$, the Fock representation can be
identified with the ring of differential polynomials
$\C[u_i^{(n)}]_{i=1,\ldots,\el;n\geq 0}$, where $u_i =
b_i(-1)$. We have proved in \cite{ff:im,ff:kdv} that the classical
limits $\Q_i, i=0,\ldots,\el$, of the screening operators act on
$\C[u_i^{(n)}]_{i=1,\ldots,\el;n\geq 0}$ and generate the pronilpotent
subalgebra $\N_+ \subset \G$. Moreover,
$\C[u_i^{(n)}]_{i=1,\ldots,\el;n\geq 0} \simeq \C[N_+/A_+]$, where
$N_+$ is the Lie group of $\N_+$ and its {\em principal} abelian
subgroup $A_+$. The left infinitesimal action of the Lie algebra
$\N_+$ on $N_+/A_+$ coincides with the one defined by the operators
$\Q_i, i=0,\ldots,\el$.

In particular, the operators $\Q_i, i=1,\ldots,\el$, generate the
action of the Lie subalgebra $\n_+ \subset \N_+$ on
$\C[u_i^{(n)}]_{i=1,\ldots,\el;n\geq 0}$, and $\WW_0(\g)$ is the ring
of invariants of this action. This ring coincides with the ring of
$\n_+$--invariants of $\C[N_+/A_+]$, and we have shown that it is also
a ring of differential polynomials
$\C[v_i^{(n)}]_{i=1,\ldots,\el;n\geq 0}$. The embedding $\C[v_i^{(n)}]
\arr \C[u_i^{(n)}]$ is called the generalized Miura transformation.

The quotient of $\C[u_i^{(n)}] \otimes \C[z,z^{-1}]$ by the total
derivatives and constants is the {\em Heisenberg-Poisson algebra}
denoted by $\F_0$.  The {\em classical $\W$--algebra} associated to
$\g$, $I_0(\g)=\W_0(\g)$ is a Poisson subalgebra of $\F_0$, which can
be identified with the quotient of $\WW_0(\g)[z,z^{-1}]$ by the total
derivatives and constants.

Next, the space $I_0(\G)$ is identified with the Lie algebra
cohomology $H^1(\N_+,\C[u_i^{(n)}])$. Using the fact that
$\C[u_i^{(n)}] \simeq \C[N_+/A_+]$, one shows that $H^*(\N_+,\pi_0)
\simeq \bigwedge^*(\ab_+^*)$, where $\ab^*_+$ is the dual space to the
Lie algebra $\ab_+$ of $A_+$. This way we we proved that $I_0(\G)
\simeq \ab_+^*$ \cite{ff:im}.

Each element of $I_0(\G)$ gives rise to a hamiltonian vector field on
$N_+/A_+$ and these vector fields the $\g$--{\em mKdV hierarchy} (on
the set of functions $u_i(t), i=1,\ldots\el$). On the other hand, the
abelian Lie algebra $\ab_-$, which is the opposite of the Lie algebra
$\ab_+$, also acts on this homogeneous space from the right by vector
fields. In \cite{ff:kdv} we proved that these two actions coincide.

Finally, these derivations preserve the differential subring
$\C[v_i^{(n)}] = \WW_0(\g)$. The corresponding equations on the set of
functions $v_i(t), i=1,\ldots\el$, form the $\g$--{\em KdV
hierarchy}.

\subsection{The present work}
In this paper we apply the same approach when the principal abelian
subgroup $A_+$ of $N_+$ is replaced by the {\em homogeneous} abelian
subgroup $H_+$. In this case we have a structure, which is very
similar to what we have in the principal case. This is illustrated by
the following table.

\bigskip

\begin{center}
\begin{tabular}{|l|l|l|}
\hline
 & Principal Case & Homogeneous Case \\
\hline
\hline
Differential polynomials & $\C[u_i^{(n)}]$ &
$\C[p_\al^{(n)},q_\al^{(n)},u_i^{(n)}]$ \\
\hline
Homogeneous space & $N_+/A_+$ & $N_+/H_+$ \\
\hline
Screening operators & $\int e^{-\phi_i} dz$ & $\int \sum_{\beta
\in \De_+} P^R_{\al_i,\beta}(q) p_\beta e^{-\phi_i} dz$ \\
\hline
Non-local equations & (affine) Toda equation & equations \eqref{new},
\eqref{nonloc1} \\
\hline
$I_0(\g)$ & classical $\W$--algebra $\W_0(\g)$ & $\C[\G^*_1]\lo$ \\
\hline
$I_\nu(\g)$ & $\W$--algebra $\W_\nu(\g)$ & $U_k(\G)\lo,
k=-h^\vee+\nu^{-1}$ \\
\hline
Embedding $I_\nu(\g) \subset \F_0^\nu$ & Miura transformation &
Wakimoto realization \\
\hline
$I_0(\G)$ & mKdV hamiltonians & mAKNS hamiltonians \\
\hline
Symmetries of the & mKdV hierarchy & mAKNS hierarchy \\

non-local equation & & \\
\hline
Local equations in $I_0(\g)$ & KdV hierarchy & AKNS hierarchy \\
\hline
\end{tabular}
\end{center}

\bigskip

Thus, the starting point of our approach is always a
Heisenberg--Poisson algebra and a pair of {\em non-local} equations,
which possess local integrals of motion (or symmetries) in this
Heisenberg--Poisson algebra. Both equations are generated by the
screening operators corresponding to the simple roots of $\g$ or
$\G$. The action of the Poisson bracket with each of these terms on
the ring of differential polynomials coincides with the action of the
corresponding generator of the pronilpotent Lie algebra $\N$ on the
homogeneous space of $N_+$. This fact allows us to interpret the local
integrals of motion as cohomologies of this Lie algebra.

In the principal case, the Poisson algebra of local integrals
of motion of the non-local equation associated to $\g$ (the $\g$--Toda
equation) is the classical $\W$--algebra $\W_0(\g)$, and the local
integrals of motion of the non-local equation corresponding to $\G$
(affine Toda equation) are the hamiltonians of the $\g$-- mKdV
hierarchy. The embedding of $\W_0(\g)$ into the Heisenberg--Poisson
algebra is the generalized Miura transformation. The equations of the
mKdV hierarchy, written in terms of $\W_0(\g)$, become the equations
of the $\g$--KdV hierarchy.

In the homogeneous case, the role of $\W_0(\g)$ is played by the
classical affine algebra, i.e., the Poisson algebra $\C[\G^*_1]\lo$ of
local functionals on a hyperplane in $\G^*$. The role of the
$\g$--mKdV hierarchy is played by a modified AKNS hierarchy for
$\g=\sw_2$ and its generalization for other $\g$, which we call the
$\g$--{\em mAKNS hierarchy}. The embedding of $\C[\G^*_1]\lo$ into the
Heisenberg--Poisson algebra is the Poisson analogue of the Wakimoto
realization. Written in terms of $\C[\G^*_1]\lo$, the $\g$--mAKNS
hierarchy becomes what we call the $\g$--{\em AKNS hierarchy}.

The next step is to quantize the classical integrals of motion. In the
principal case, this was done in \cite{ff:im}. Here we prove the
existence of the quantum integrals of motion of the non-local equation
corresponding to $\g$ in the homogeneous case. This gives us an
embedding of $\G$ into a Heisenberg algebra, analogous of the
embedding of the $\W$--algebra $\W_\nu(\g)$ into $\F^\nu_0$. This free
field realization of $\G$ is nothing but the Wakimoto realization,
which has been previously defined in \cite{Wak} for $\g=\sw_2$ and in
\cite{ff:usp} for general $\g$. Thus we obtain a new proof of the
Wakimoto realization of affine algebras. The original construction of
the Wakimoto realization has been geometric in nature, see
\cite{ff:usp}. In this paper we give another perspective on the
Wakimoto realization, which comes from the theory of integrable
hierarchies of soliton equations. This confirms to the general
philosophy of ``free field realization'' outlined in \cite{Fr}.

Finally, we conjecture that all classical integrals of motion of the
non-local equation corresponding to $\G$ can be quantized.

The main results of this paper have been obtained during our visit to
Kyoto University in the Summer of 1993, and the present paper is an
expanded version of a draft written in the Fall of 1994. In the mean
time, some aspects of the $\sw_2$ case of our program have been
considered in \cite{ABF}. We also became aware of the interesting
papers \cite{To}, in which integrable hierarchies similar to ours have
been considered, from a different point of view.

The program outlined above has been further developed in
\cite{EF1,EF2,Five}. In \cite{BF} the results of \cite{ff:kdv} have
been interpreted geometrically in terms of moduli spaces of bundles,
thus allowing a uniform treatment of the integrable hierarchies
associated to arbitrary maximal abelian (or Heisenberg) subalgebras of
$\G$.

\section{Wakimoto realization}    \label{wreal}

Let $\g$ be a simple Lie algebra with the Cartan decomposition $\g = \n_+
\oplus \h \oplus \n_-$, where $\n_+$ and $\n_-$ are the upper and lower
nilpotent subalgebras, respectively, and $\h$ is its Cartan subalgebra.

Let $\G$ be the affine algebra, corresponding to $\g$: the universal
central extension of the loop algebra $\g[t,t^{-1}]$.

Introduce the Heisenberg algebra $\HH(\g)$, which has generators $a_\al(n),
a^*_\al(n), \al \in \De_+, n \in \Z$, where $\De_+$ is the set of positive
roots of $\g$, and $\one$. The commutation relations read:
$$[a_\al(n),a_\beta^*(m)] = \delta_{\al,\beta} \delta_{n,-m} \one, \quad
\quad [a_\al(n),a_\beta(m)] = 0, \quad \quad [a_\al^*(n),a_\beta^*(m)] =
0,$$ and $\one$ commutes with everything.

For $\nu \neq 0$, denote by $M$ the Fock representation of $\HH(\g)$,
which is generated from the vacuum vector $v$ satisfying $$a_\al(n) v = 0,
\quad n\geq 0; \quad \quad a^*_\al(n) v = 0, \quad n>0,$$ and on which the
central element $\one$ acts as the identity.

The space $M$ is a vertex operator algebra (VOA) \cite{B,FLM}. We will give
an explicit formula for the currents $Y(\cdot,z)$ defined by this VOA.

Monomials $a_{\al_1}(m_1) \ldots a_{\al_k}(m_k) a_{\beta_1}^*(n_1) \ldots
a^*_{\beta_l}(n_l) v, m_p < 0, n_q\leq 0$, form a linear basis in
$M$. The series $Y(\cdot,z)$ associated to this monomial is given by
$$C \, \,:\pa_z^{-m_1-1} a_{\al_1}(z) \ldots \pa_z^{-m_k-1} a_{\al_k}(z)
\pa_z^{-n_1} a_{\beta_1}^*(z) \ldots \pa_z^{-n_l} a_{\beta_l}^*(z):,$$
where $C = [(-m_1-1)! \ldots (-m_1-1)! (-n_1)! \ldots (-n_l)!]^{-1}$,
columns stand for normal ordering, and
\begin{equation}    \label{betgam}
a_\al(z) = \sum_{n \in \Z} a_\al(n) z^{-n-1}, \quad \quad a_\beta^*(z) =
\sum_{n\in\Z} a_\beta^*(n) z^{-n}.
\end{equation}

The Fourier coefficients of currents of this VOA form a Lie algebra,
which lies in a certain completion of the universal enveloping algebra
$U(\HH(\g))$ factored by the ideal generated by $(\one - 1)$,
respectively. Following \cite{ff:gd}, we call this Lie algebra the
{\em local completion} of the universal enveloping algebra and denote
it by $U(\HH(\g))_{\on{loc}}$.

Denote by $L\n_+$ the Lie subalgebra $\n_+[t,t^{-1}]$ of $\G$.
Consider the Lie group $LN_+$ of $L\n_+$. It consists of polynomial maps
from the circle to the nilpotent subgroup $\NN_+$ of $G$.

The Lie group $\NN_+$ is isomorphic to its Lie algebra $\n_+$ via the
exponential map. This allows us to introduce coordinates $x_\al, \al \in
\De_+$, on $\NN_+$, such that $x_\al$ has weight $\al$ with respect to the
natural action of the Cartan subgroup $H \subset G$, which is the Lie group
of $\h$. Denote the vector field corresponding to the left (respectively,
right) action of a generator $e_\al$ of $\n_+$ on $\NN_+$ by $e_\al^L$
(respectively, $e_\al^R$). In coordinates $x_\beta$ they are given by:
\begin{equation}    \label{LR}
e_\al^L = \sum_{\beta \in \De_+} P^L_{\al,\beta} \frac{\pa}{\pa x_\beta},
\quad \quad
e_\al^R = \sum_{\beta \in \De_+} P^R_{\al,\beta} \frac{\pa}{\pa x_\beta},
\end{equation}
where $P^L_{\al,\beta}$ and $P^R_{\al,\beta}$ are certain polynomials in
$x_\gamma, \gamma \in \De_+$ of degree $\beta-\al$.

Coordinates $x_\al, \al \in \De_+$, on the group $\NN_+$ give us coordinates
$x_\al(n), \al \in \De_+, n \in \Z$, on the group $L\NN_+$.  Denote
$a^*_\al(n) = x_\al(-n), a_\al(n) = \pa/\pa x_\al(n)$. These operators
generate the Heisenberg algebra $\HH(\g)$. We have two commuting
infinitesimal actions of $L\n_+$ on $L\NN_+$ by vector fields: left and
right.

Explicitly, we have
\begin{equation}    \label{L}
e^L_\al(z) = \sum_{n\in\Z} e^L_\al(n) z^{-n-1} = \sum_{\beta \in
\De_+} P^L_{\al,\beta}(z) a_\beta(z),
\end{equation}
\begin{equation}    \label{R}
e^R_\al(z) = \sum_{n\in\Z} e^R_\al(n) z^{-n-1} = \sum_{\beta
\in \De_+} P^R_{\al,\beta}(z) a_\beta(z),
\end{equation}
where $P^L_{\al,\beta}(z)$ and $P^R_{\al,\beta}(z)$ are obtained from the
polynomials $P^L_{\al,\beta}$ and $P^R_{\al,\beta}$, respectively, by
replacing $x_\gamma, \gamma \in \De_+$, by $a_\gamma^*(z)$, and we use the
notation $e_\al(n) = e_\al \otimes t^n$.

These formulas define two embeddings $L\n_+ \arr U(\HH(\g))_{\on{loc}}$ and
hence two commuting actions of the Lie algebra $L\n_+$ on the space $M$.

\begin{rem}
We can define these actions in a coordinate independent way.

Introduce the Lie algebra $D(\n_+)$. It has generators $y^R(n), y \in \n_+,
n \in \Z$, and $P(n), P \in \C[\NN_+], n \in \Z$, where $\C[\NN_+]$ stands for
the ring of regular functions on $\NN_+$. If we choose coordinates $x_\al$ on
$\NN_+$, then $\C[\NN_+] \simeq \C[x_\al]_{\al\in\De_+}$. There are the
following relations: $$[y_1^R(n_1),y_2^R(n_2)] = [y_1,y_2]^R(n_1+n_2),
\quad \quad [y^R(n),P(m)] = [y \cdot P](n+m),$$ where $y \cdot P$ denotes
the action of $y \in \n_+$ on $P \in
\C[\NN_+]$ by vector field from the right, and $$[P(n),Q(m)] = 0, \quad \quad
P(n) = \sum_{n_1+n_2=n} P_1(n_1) P_2(n_2),$$ if $P = P_1 P_2$ in
$\C[\NN_+]$. If we choose coordinates on $\NN_+$, the Lie algebra $D(\n_+)$
becomes isomorphic to the Heisenberg algebra $\HH(\g)$.

The linear span of $y^R(n), y \in \n_+, n \geq 0$, and $P(n), P \in
\C[\NN_+], n>0$, is a Lie subalgebra $D_+(\n_+)$ of $D(\n_+)$. The
module $M$ over $D(\n_+)$ can be defined as the module induced from
the trivial one-dimensional representation $\C v$ of $D_+(\n_+)$.

The correlation functions, i.e. matrix elements of the currents $$y^R(z) =
\sum_{n\in\Z} y^R(n) z^{-n-1}, \quad \quad P(z) = \sum_{n\in\Z} P(n)
z^{-n}$$ can also be expressed in terms of action of $\n_+$ on $\C[\NN_+]$ in
a coordinate independent way. These correlation functions uniquely
determine vertex algebra structure on $M$.

One obtains, e.g., the following formula (compare with \cite{vs,aty}):
$$\left\lef v^* , \prod_{s=1}^m e^R_{i_s}(w_s) \prod_{j=1}^N P_j(z_j) v
\right\ri = \sum_{p=(I^1,\ldots,I^N)} \prod_{j=1}^N \frac{\jmath(
e_{i^j_{a_j}}^R \ldots e_{i^j_1}^R P_{X_j} )}{(w_{i^j_1}-w_{i^j_2})\ldots
(w_{i^j_{a_j}}-z_j)},$$ where the summation is taken over all {\em ordered}
partitions $I^1 \cup I^2 \cup \ldots \cup I^N$ of the set $\{
i_1,\ldots,i_m \}$, where $I^j = \{ i^j_1,i^j_2,\ldots,i^j_{a_j} \}$, and
$\jmath$ is the augmentation homomorphism $\C[\NN_+] \arr \C$.

We note that the construction described above assigns a vertex algebra
to an arbitrary affine algebraic group in place of $\NN_+$.\qed
\end{rem}
\smallskip

Now put $W^\nu_\la = M \otimes \pi^\nu_\la$, where $\pi^\nu_\la$ was
defined in the Introduction. The space $W^\nu_0$ is a VOA, the tensor
product of the VOAs $M$ and $\pi^\nu_0$. We can extend the action of
$\pa$ to $W^\nu_\la[z,z^{-1}]$ by the formula $\pa \otimes 1 + 1
\otimes \pa_z$.  Denote by $\w^\nu_\la$ the quotient of $W^\nu_\la
\otimes \C[z,z^{-1}]$ by the total derivatives and constants, if
$\la=0$. The space $\w^\nu_\la$ is a Lie algebra, which is isomorphic
to $U_\nu(\HH(\g) \oplus \hh)_{\on{loc}}$. Introduce operators
$D^\nu_i(n): W^\nu_\la \arr W^\nu_{\la-\al_i}, i=1,\ldots,\el$, by the
formula
\begin{equation}    \label{di}
D^\nu_i(z) = \sum_{n\in\Z} D^\nu_i(n) z^{-n-\nu(\al_i,\la)} =
e^R_{i}(z) V^\nu_{-\al_i}(z),
\end{equation}
where $V^\nu_{-\al_i}(z)$ is given by formula \eqref{vertex}.

Put $$G^\nu_i = D^\nu_i(1) = \sum_{n\in\Z} e^R_{i}(n)
V^\nu_{-\al_i}(-n), \quad \quad G^\nu_i: W^\nu_0 \arr
W^\nu_{-\al_i}.$$ This operator is called the $i$th screening
operator. It commutes with the action of $\pa$. We denote by
$\ol{G}^\nu_i$ the induced operator $\w^\nu_0 \arr \w^\nu_{-\al_i},
i=1,\ldots,\el$.

Denote $$K_\nu(\g) = \bigcap_{i=1}^\el \on{Ker}_{W^\nu_0} G^\nu_i,
\quad \quad J_\nu(\g) = \bigcap_{i=1}^\el \on{Ker}_{\w^\nu_0}
\ol{G}^\nu_i.$$ According to Lemma 4.2.8 from \cite{ff:im},
$K_\nu(\g)$ is a vertex algebra, and $J_\nu(\g)$ is a Lie algebra.

Let $V_k, k \in \C$, be the VOA of the affine algebra $\G$. Recall
that as a vector space $$V_k = U(\G) \otimes_{U(\g[t] \oplus \C K)}
\C_k,$$ where $\C_k$ stands for the trivial one-dimensional
representation of the Lie subalgebra $\g[t]$ of $\G$, on which $K$
acts by multiplication by $k$. Its $\Z$--grading is inherited from the
standard $\Z$--grading on $\G$, such that $\deg A(n) = -n, \deg K =
0$, cf.  \cite{ff:gd,fzhu}.

The Fourier coefficients of currents of the VOA $V_k$ form a Lie algebra
$U_k(\G)_{\on{loc}}$, the local completion of the universal enveloping
algebra at level $k$, $U_k(\G) = U(\G)/(K-k)U(\G)$.

\begin{thm}    \label{generic}
For generic $\nu$, the vertex algebra $K_\nu(\g)$ is isomorphic to the
VOA $V_k$ of the affine algebra $\G$, and $J_\nu(\g)
\simeq U_k(\G)_{\on{loc}}$, where $k = - h^\vee + \nu^{-1}$, $h^\vee$ being
the dual Coxeter number of $\g$.
\end{thm}

The proof of this theorem is analogous to the proof of Theorem 4.5.9 from
\cite{ff:im}. We construct a family of complexes $C^*_\nu(\g)$ depending on
the parameter $\nu$, whose $0$th cohomology is $K_\nu(\g)$ and which has a
well-defined {\em classical limit}, when $\nu \arr 0$. We then show that in
this limit all higher cohomologies of the complex vanish and the $0$th
cohomology can be identified with the limit of the VOA $V_k$ when $k
\arr \infty$. This will allow us to compute the cohomology of the complex
$C_\nu(\g)$ for generic $\nu$ and identify $K_\nu(\g)$ with the VOA $V_k$.

\medskip

\noindent{\em Example}. We give here an explicit realization of the
kernel of the screening operator $G_1$ in the case when $\g=\sw_2$
(for generic $\nu$). In the following formulas we suppress the index
$1$. Let $\{ e,h,f \}$ be the standard basis of $\sw_2$. Set
\begin{align}    \notag
e(z) &= a(z), \quad \quad h(z) = -2 :a(z) a^*(z): + \frac{1}{\nu}
b(z),\\    \label{f(z)}
f(z) &= - :a(z) a^*(z) a^*(z): + (-2+\nu^{-1}) \pa_z a^*(z) +
\frac{1}{\nu} b(z) a^*(z).
\end{align}
These formulas first appeared in \cite{Wak}. The Fourier coefficients
of the above series satisfy the relations in $\su$ with $k = -2 +
\nu^{-1}$. Moreover, the generating vector of $W^\nu_0$ is annihilated
by all non-negative Fourier coefficients. Therefore we obtain a
homomorphism $V_k \arr W^\nu_0$. It is known that $V_k$ is irreducible
for generic $k$. Hence for generic $\nu$ the above homomorphism is
injective. On the other hand, it is shown in \cite{ff:lmp} that the
generating series above commute with the screening operator $$G_1^\nu
= - \sum_{n \in \Z} a(n) V^\nu_{-\al_1}(-n): W^\nu_0 \arr
W^\nu_{-2}.$$ Therefore the image of $V_k$ in $W^\nu_0$ lies in the
kernel of $G_1$. The calculation of the characters given below proves
that $V_k$ is equal to the kernel of $G_1$.

In \cite{ff:usp,ff:lmp} we generalized the above construction to the
case of an arbitrary simple Lie algebra $\g$. The commutativity with
the screening operators was proved in \cite{Fr:th,FFR}. This gives us
a proof of \thmref{generic}. In what follows, we will give an
alternative proof of this theorem.

\section{The complex}    \label{complex}
The $j$th group $C^j_\nu(\g)$ of our complex is $$C^j_\nu(\g) =
\oplus_{l(s)=j} W^\nu_{s(\rho)-\rho},$$ where $s$ belongs to the Weyl
group, $l(s)$ is its length, and $\rho \in \h^*$ is the half-sum of
positive roots of $\g$.

The construction of the differentials of the complex $C^*_\nu(\g)$ follows
closely the construction of the differentials of the complex
$F^*_\nu(\g)$ from \cite{ff:im}, Sect.~4.5 (where $\beta^2$ played the role
of $\nu$).

Let $p = (p_1,\ldots,p_m)$ be a permutation of the set $(1,2,\ldots,m)$. We
define a contour of integration $C_p$ in the space $(\C^\times)^m$ with the
coordinates $z_1,\ldots,z_m$ as the product of one-dimensional contours
along each of the coordinates, going counterclockwise around the origin
starting and ending at the point $z_i = 1$, and such that $|z_{p_1}| >
|z_{p_2}| > \ldots > |z_{p_m}|$ whenever $z_i \neq 1$.

Denote by ${\mathbf i} = (i_1,\ldots,i_m)$ a sequence of numbers from $1$ to
$l$, such that $i_1 \leq i_2 \leq \ldots \leq i_m$. We can apply a
permutation $p$ to this sequence to obtain another sequence $p(\is) =
(i_{p(1)},\ldots,i_{p(m)})$. Put $\gamma = \sum_{j=1}^m \al_{i_j}$.  Let us
define an operator $D_{p(\is)}^\nu$ from $W^\nu_\la$ to
$W^\nu_{\la-\gamma}$ as the integral $$\int_{C_p} dz_1 \ldots dz_m \, \,
D_{i_{p(1)}}(z_1) \ldots D_{i_{p(m)}}(z_m) =$$ $$\int_{C_p} dz_1 \ldots
dz_m \prod_{1\leq k<l\leq m} (z_k - z_l)^{\nu(\al_{i_k},\al_{i_l})}
\prod_{1\leq k\leq m} z_k^{\nu(\la,\al_{i_k})}
:\V^\nu_{-\al_{i_1}}(z_1) \ldots \V^\nu_{-\al_{i_m}}(z_m): \cdot$$ $$\cdot {\mc
E}_{i_1,\ldots,i_m}(z_1,\ldots,z_m),$$ where $$\V^\nu_{-\al_i}(z) =
\sum_{n\in\Z} V^\nu_{-\al_i}(n) z^{-n},$$ and $${\mc
E}_{i_1,\ldots,i_m}(z_1,\ldots,z_m) = e^R_{{i_1}}(z_1) \ldots
e^R_{{i_m}}(z_m).$$

The latter can be rewritten using Wick's formula as a linear combination of
normally ordered products of currents $a_\al(z)$ and $a^*_\al(z), \al \in
\De_+$, multiplied by rational functions in $z_1,\ldots,z_m$, which have
poles only on the diagonals $z_i=z_j$.

$D^\nu_{p({\mathbf i})}$ is a linear operator from $W^\nu_\la$ to the
completion $\widehat{W}^\nu_{\la+\gamma}$ of $W^\nu_{\la+\gamma}$. Note
that this operator is uniquely defined by the sequence ${\mathbf j} =
p(\is)$, and so we can denote it by $D_{{\mathbf j}}^\nu$.

We choose the branch of a power function appearing in the integral, which
takes real values for real $z_i$'s ordered so that $z_{j_1} > z_{j_2} >
\ldots > z_{j_m}$. Thus, $C_p$ should be viewed as an element of the group
of relative $m$-chains in $(\C^\times)^m$ modulo the diagonals, with values
in the one-dimensional local system $\xi_{\js}$, which is defined by the
multi-valued function $$\prod_{1\leq k<l\leq m} (z_k -
z_l)^{-\nu(\al_{i_k},\al_{i_l})} \prod_{1\leq k\leq m}
z_k^{-\nu(\la,\al_{i_k})}.$$

Our integral is well-defined for generic values of $\nu$ over any such
relative chain. Indeed, the integral $$\int_{C} dz_1 \ldots dz_m
\prod_{1\leq k<l\leq m} (z_k - z_l)^{\mu_{kl}}
\prod_{1\leq k\leq m} z_k^{\nu_k}$$ over such a chain $C$ converges in the
region $\mbox{Re} \mu_{kl} \geq 0$, and can be analytically continued to
other values of $\mu_{kl}$, which do not lie on hyperplanes
\begin{equation}    \label{hyperplane}
\sum_{k,\el \in S,k<l} \mu_{kl} = -s, \quad \quad s \in \Z, s \geq (\# S)-1,
\end{equation}
where $S$ is a subset of the set $\{ 1,2,\ldots,m \}$,
cf. \cite{varchenko}, Theorem (10.7.7), for details.

For generic $\nu$ the exponents in our integral do not lie on those
hyperplanes. Indeed, an expression of the form $$\sum_{1 \leq k<l \leq r}
\nu (\al_{j_k},\al_{j_l})$$ can take integral value for generic $\nu$,
if and only if
\begin{equation}    \label{zero}
\sum_{1 \leq k<l \leq r} (\al_{j_k},\al_{j_l}) = 0.
\end{equation}
If this is the case, for our integral to converge, the rational functions
appearing in ${\mc E}_{j_1,\ldots,j_r}(z_1,\ldots,z_r)$ should not have
poles on the diagonals of combined order $r-1$ or more.

In the case when $\al_{j_1} = \al_{j_2} = \ldots = \al_{j_{r-1}} =
\al_i$, and $\al_{j_r} = \al_j$, this fact follows from the Serre
relations: $$(\on{ad} e^R_i)^{-a_{ij}+1} \cdot e_j(m) = 0,$$ which the
operators $e^R_i$ satisfy. Indeed, in order to satisfy
\eqref{zero}, we should put $r=-2(\al_i,\al_j)/(\al_i,\al_i)+2$. But the
Serre relations imply that the coefficients of ${\mc
E}_{i,\ldots,i,j}(z_1,\ldots,$ $z_r)$ are rational functions in
$z_1,\ldots,z_r$, whose combined poles have combined order less than or
equal to $-a_{ij} = -2(\al_i,\al_j)/(\al_i,\al_i) < r-1$. Therefore, the
equation \eqref{hyperplane} can not hold in this case.

In general we have to show that if \eqref{zero} holds then any commutator
of the form $[e_{m_1},[e_{m_2},[...,e_{m_r}]\ldots]]$, where
$(m_1,\ldots,m_r)$ is a permutation of the set $(j_1,\ldots,j_r)$,
vanishes. This is indeed the case, since $\gamma = \sum_{k=1}^r
\al_{j_r}, r>1$, can not be a root of $\g$, if \eqref{zero} holds.

The proof, which works for an arbitrary Kac-Moody algebra with a
symmetrizable Cartan matrix, was communicated to us by V. Kac
(cf. \cite{K}): the equation \eqref{zero} can be rewritten as
$2(\rho,\gamma) = (\gamma,\gamma)$. This equation can not hold for
imaginary roots, because then $(\gamma,\gamma)\leq 0$ and
$(\rho,\gamma)>0$. If $\gamma$ is real, i.e. $\gamma = w \cdot \al_i$
for a simple root $\al_i$ and an element $w$ of the Weyl group of
$\g$, then we obtain $2(\rho,w \cdot \al_i) = (\al_i,\al_i)$. This is
true for $w=1$, but each reflection from the Weyl group increases the
left hand side because $w \cdot \al_i > 0$ by assumption while the
right hand side remains the same. Hence this equality can not hold for
$w \neq 1$.

This shows that our integrals are well-defined.

We can interpret the operator $D^\nu_{{\mathbf j}}$ as a composition
operator $G^\nu_{j_1} \ldots G^\nu_{j_m}$, cf. \cite{ff:im},
Sect.~4.5.3. The operators $G^\nu_j$ satisfy the $q$--Serre relations,
where $q = \exp (\pi i \nu)$, in the following sense \cite{bmp} (see
also \cite{bmp1,ff:im}).

Consider a free algebra $A$ with generators $g_i, i=1,\ldots,\el$. We can
assign to each monomial $g_{j_1} \ldots g_{j_m}$ the contour $C_{\js}$,
where $\js = (j_1,\ldots,j_m)$, and hence the operator $D_{\js}^\nu$. This
gives us a map $\Delta$ from $A$ to the space of linear combinations of
such contours. Given such a linear combination $C$, we define $D_C^\nu$ as
the linear combination of the corresponding operators $D_{C_{\js}}^\nu$.

Consider the two-sided ideal $S_q$ in $A$, which is generated by the
$q$-Serre relations $(\mbox{ad} g_i)_q^{-a_{ij}+1} \cdot g_j,
i\neq j$, where $q=\exp (\pi i \nu)$.

\begin{lem}    \label{qserre}
If $C$ belongs to $\Delta(S_q)$, then $D_C^\nu = 0$.
\end{lem}

The proof is given in \cite{bmp1} (note however that in that paper the
question of convergence of integrals was not addressed). It is based
on rewriting the integrals over the contours $C_{\js}$ as integrals
over other contours, where all variables are on the unit circle with
some ordering of their arguments.

\lemref{qserre} means that the operators $G_i^\nu$ ``satisfy'' the
$q$--Serre relations of $\g$. Thus, we obtain a well-defined map, which
assigns to each element $P$ of the algebra $U_q(\n_+) \simeq A/S_q$ the
operator $D_P^\nu$.

\begin{lem}    \label{sing}
Let $P \in U_q(\n_-)$ be such that $P \cdot {\mathbf 1}_\la$ is a singular
vector of weight $\la+\gamma$ in the Verma module $M_\la^q$ of highest
weight $\la$ over $U_q(\g)$. Then the operator $D_P^\nu$ is a homogeneous
linear operator $W^\nu_\la \arr W^\nu_{\la+\gamma}$, which commutes with
the action of $\pa$.
\end{lem}

The proof is given in \cite{ff:im}, Sect.~4.5.6.

We are ready now to define the differentials $\delta^j_\nu: C^j_\nu(\g)
\arr C^{j+1}_\nu(\g)$ of the quantum complex $C^*_\nu(\g)$. Recall that for
any pair $s, s'$ of elements of the Weyl group of $\g$ there exists a
singular vector $P_{s',s}^q
\cdot v_{s(\rho)-\rho}^q$ of weight $s'(\rho)-\rho$ in the Verma module
$M^q_{s(\rho)-\rho}$, cf. \cite{ff:im}, Sect.~4.4.5. We put:
\begin{equation}    \label{differential}
\delta^j_\nu = \sum_{l(s)=j,l(s')=j+1,s\prec s'} \ep_{s',s} \cdot
D_{P_{s',s}^q}^\nu,
\end{equation}
where $q = \exp (\pi i \nu)$, and $\ep_{s',s} = \pm 1$ are signs chosen in
a special way. By \lemref{sing}, the differentials $\delta^j_\nu$ are
well-defined homogeneous linear operators. From the nilpotency of the
differential of the quantum BGG resolution, cf. \cite{ff:im}, Sect.~4.4.6, and
\lemref{qserre} we derive that these differentials are nilpotent:
$\delta^{j+1}_\nu \delta^j_\nu = 0$.

Thus, we have constructed a family of complexes $C^*_\nu(\g)$. We have:
$C^0_\nu(\g) = W^\nu_0, C^1_\nu(\g) = \oplus_{i=1}^\el W^\nu_{-\al_i}$,
and $\delta^0: C^0_\nu(\g) \arr C^1_\nu(\g)$ is given by the sum of the
operators $G^\nu_i: W^\nu_0 \arr W^\nu_{-\al_i}$. Therefore the $0$th
cohomology of the complex $C^*_\nu(\g)$ is nothing but the VOA $K_\nu(\g)$.

Since the differentials of the complex $C^*_\nu(\g)$ commute with the
action of $\pa$, we can form the double complex $$\C \larr C^*_\nu(\g)
\stackrel{\pm \pa}{\larr} C^*_\nu(\g) \larr \C.$$ The $0$th cohomology of
the total complex $\cc^*_\nu(\g)$ of this double complex is the
space $J_\nu(\g)$.

In order to compute the cohomologies of the complexes $C^*_\nu(\g)$ and
$\cc^*_\nu(\g)$, we will study their {\em classical limit} $\nu
\arr 0$.

\begin{rem}    \label{ser}
Let $\g$ be an arbitrary symmetrizable Kac-Moody algebra, and
$\al_1,\ldots,\al_l$ be the set of simple roots of $\g$. Consider a
set of homogeneous linear operators $X_i(n), i=1,\ldots,\el; n \in
\Z$, acting on a $\Z$--graded linear space $M$, so that $\deg X_i(n) =
-n$. Define the series $X_i(z) V^\nu_{-\al_i}(z), i=1,\ldots,\el,$ of
linear operators acting from $M \otimes \pi^\nu_\gamma$ to $M \otimes
\pi^\nu_{\gamma-\al_i}$. We can define integrals of products of these
operators in the same way as above. Our proof of convergence above
applies in the general case as well, and it shows that these integrals
converge and \lemref{qserre} holds, if and only if $X_i(n)$'s satisfy
the Serre relations of $\g$:
$$[X_i(n_1),[X_i(n_2),...,[X_i(n_{-a_{ij}+1}),X_j(m)]...]] = 0, \quad \quad
n_i,m \in \Z.$$ In other words, the operators $\int X_i(z)
V^\nu_{-\al_i}(z) dz$ ``satisfy'' the $q$--Serre relations, if and only if
the operators $X_i(n)$ satisfy the Serre relations.\qed
\end{rem}

\section{Classical limit}    \label{class1}
Introduce new operators $a'_\al(n) = \nu a_\al(n)$.  We have the
commutation relations $$[a'_\al(n),a^*_\beta(m)] = \nu \delta_{\al,\beta}
\delta_{n,-m}$$ for the operators $a'_\al(n)$ and $a^*(m)$ acting on $M$.

Consider a linear basis in $M$, which consists of monomials in $a'_\al(n),
\al \in \De_+, n<0$, and $a^*_\al(n), \al \in \De_+, n\leq 0$, applied to
the vacuum vector $v$. As a basis in $\pi^\nu_\la$ we take monomials in
$b_i(n), i=1,\ldots,\el, n<0$, applied to the vacuum vector $v_\la$. Tensor
products of elements of these bases form a basis in $W^\nu_\la$. Using these
bases, we can identify the spaces $W^\nu_\la$ with different values of
$\nu$, so we can omit the superscript $\nu$ and write $W_\la$. The linear
operators defined above for different $\nu$ should be considered as
operators explicitly depending on the parameter $\nu$, acting spaces
$W_\la$, which do not depend on $\nu$.

We can also identify the spaces $\w^\nu_\la$ with different $\nu$ and
write $\w_\la$. However the Lie algebra structure on $\w_0$
will be $\nu$--dependent.

Consider the operator $G_i^\nu: W_0 \arr W_{-\al_i}$. It can be
expanded in powers of $\nu$: $G_i^\nu = G_i^0 + \nu(\ldots)$. We will
call $G_i^0$ the classical screening operators. We can consider
$G_i^0$ as an operator $W_\la \arr W_{\la-\al_i}$ for any $\la$.

Put $G_i = T_{\al_i} G_i^0: W_0 \arr W_0$, where $T_{\al_i}: W_{-\al_i}
\arr W_0$ was defined after formula \eqref{vertex}.

\begin{lem}    \label{nilp1}
The operators $G_i$ generate an action of the nilpotent subalgebra $\n_+$
of $\g$ on $W_0$.
\end{lem}

The Lemma follows from \lemref{qserre} in the limit $q \arr 1$. A different
proof will be given below.

In the same way as in \cite{ff:im} we can show that $D_{P_{s's}^q}^\nu =
P_{s',s}(G) + \nu(\ldots)$, where $P_{s',s} \in U(\n_-)$ is the expression
for the singular vector at $q=1$, and $P_{s',s}(G)$ is the operator
$W_{s(\rho)-\rho} \arr W_{s'(\rho)-\rho}$ obtained by inserting $G_i^0$
instead of the $e_i$ in $P_{s',s}$ for all $i=1,\ldots,\el$.

Thus, we have a well-defined limit of the complex $C^*_\nu(\g)$ as $\nu
\arr 0$. As a linear space it does not depend on $\nu$: $C^j_0(\g) =
\oplus_{l(s)=j} W_{s(\rho)-\rho}$, and the differential is the $\nu \arr 0$
limit of the differential \eqref{differential}: $$\delta^j_\nu =
\sum_{l(s)=j,l(s')=j+1,s\prec s'} \ep_{s',s} \cdot P_{s',s}(G).$$

In the same way as in \cite{ff:im} we obtain the following result.

\begin{prop} The cohomology of the complex $C^*_0(\g)$ is isomorphic to the
cohomology of $\n_+$ with coefficients in $W_0$, $H^*(\n_+,W_0)$.
\end{prop}

The action of $\n_+$ on $W_0$ has geometric origin. The space $W_0$ can be
considered as the algebra of regular functions on an infinite-dimensional
linear space ${\mc U}$ with coordinates $a'_\al(n), \al \in \De_+, n<0,
a^*_\al(n), \al \in \De_+, n\leq 0$, and $b_i(n), i=1,\ldots,\el, n<0$. The
Lie algebra $\n_+$ acts on this space by vector fields. This is the
infinitesimal action of $\n_+$, corresponding to an action of the Lie group
$\NN_+$ on ${\mc U}$. We will see later that this action of $\NN_+$ is free.
Therefore the $\n_+$--module $W_0$ is co-free, i.e. dual to a free module.
Hence $H^i(\n_+,W_0) = 0, i\neq 0$, and $H^0(\n_+,W_0)$ is the algebra
$\C[{\mc U}]^{\NN_+}$ of $\NN_+$--invariant functions on ${\mc U}$. The
latter is a polynomial algebra of infinitely many variables.

We can find the degrees of the generators of the algebra $\C[{\mc
U}]^{\NN_+}$ by computing its character, which coincides with the Euler
character of the complex $C^*_0(\g)$. We obtain:
\begin{equation}    \label{char}
\on{ch} \C[{\mc U}]^{\NN_+} = \prod_{n=1}^\infty (1-q^n)^{-l} \prod_{\al \in
\De} (1-q^n u^{\al})^{-1}.
\end{equation}
This gives us the following result.

\begin{thm}    \label{zeroth}
The $0$th cohomology $K_0(\g)$ of the complex $C_0^*(\g)$ is
isomorphic to a graded polynomial algebra of infinitely many variables,
whose character is given by formula \eqref{char}. All higher cohomologies
of the complex $C^*_0(\g)$ vanish.

The $0$th cohomology $J_0(\g)$ of the complex $\cc_0^*(\g)$ is
isomorphic to the quotient of $K_0(\g)[t,t^{-1}]$ by the
subspace of total derivatives and constants. All higher cohomologies of the
complex $\cc_0^*(\g)$ vanish.
\end{thm}

\begin{cor}    \label{gen}
For generic $\nu$ all higher cohomologies of the complex $C^*_\nu(\g)$
vanish, and the character of the $0$th cohomology $K_\nu(\g)$ is given by
formula \eqref{char}.

The $0$th cohomology $J_\nu(\g)$ of the complex $\cc_0^*(\g)$ is
isomorphic to the quotient of $K_\nu(\g)[z,z^{-1}]$ by the
subspace of total derivatives and constants. All higher cohomologies of the
complex $\cc_\nu^*(\g)$ vanish.
\end{cor}

In the next section we will identify the $0$th cohomology of the complex
$C^*_0(\g)$ with the classical limit of the VOA $V_k$ of $\G$. This will
complete the proof of \thmref{generic}.

\begin{rem}    \label{limit}
\thmref{zeroth} can be proved in a simpler way by considering a
different classical limit of the differentials of the complex. Namely,
we can take the linear basis in $W_\la^\nu$, which consists of
monomials in $a_\al(n), a^*_\al(n)$, and $b'(n) = \nu^{-\frac{1}{2}}
b(n)$. Then in the limit $\nu \arr 0$ the vertex operator
$V_\gamma^\nu(z) \arr \on{Id}$, and so $G_i^\nu \arr
e_i^R(0)$. Therefore the cohomology of the complex $C^*_\nu(\g)$ in
this limit become $H^*(\n_+,\pi_0)$, where now $\n_+$ acts on $\pi_0$
via the operators $e_\al^R(0), \al \in \De_+$. This action is co-free
by definition of the operators $e_\al^R$. Hence we obtain
\thmref{zeroth}.

This proof is simpler than the proof given above, but is not clear how to
identify $V_k$ with the $0$th cohomology of $C^*_\nu(\g)$ for $\nu \neq
0$. For that it is more convenient to use the other classical limit
introduced above.\qed
\end{rem}

\section{Hamiltonian structure}
Consider the loop spaces $L\n_+ = \n_+[t,t^{-1}]$ and $L\h = \h
\otimes \C[t,t^{-1}]$. We have the inner product $\lef\cdot,\cdot\ri$
on the direct sum $L\n_+ \oplus (L\n_+)^*$, induced by the pairing
$L\n_+ \times (L\n_+)^* \arr \C$, and the inner product on $L\h$,
which are the restrictions of the invariant inner product on
$\g[t,t^{-1}]$:
$$\lef u(t),v(t) \ri = \int (u(t) , dv(t)),$$ where $(\cdot,\cdot)$ is
the invariant inner product on $\g$ normalized so that the square of
the maximal root equals $2$.

Denote by $p_\al, \al \in \De_+$ coordinates on $\n_+$, by $q_\al,
\al \in \De_+$, the dual coordinates on $\n_+^*$, and by $u_i,
i=1,\ldots,\el$, be coordinates on $L\h$. Let $\w_0$ be the space of
{\em local functionals} on $L\n_+ \oplus (L\n_+)^* \oplus L\h$. It consists
of functionals of the form $$\int
P(p^{(n)}_\al(t),q^{(n)}_\al(t),u^{(n)}_i(t);t) dt,$$ where $P$ is a
polynomial in $(p_\al(t),q_\al(t),u_i(t)) \in L\n_+
\oplus (L\n_+)^* \oplus L\h$ and their derivatives, and $t$. We use
notation $f^{(n)}(t) = \pa^n f(t)$.

Local functionals can be considered as infinite sums of monomials in the
Fourier coefficients $p_\al(n) = \int p_\al(t) t^n dt, q_\al(n) = \int
q_\al(t) t^{n-1} dt, u_i(n) = \int u_i(t) t^n dt$ of the polynomials
$p_\al(t), q_\al(t)$, and $u_i(t)$, cf. \cite{ff:im}, Sect.~2.1.

Put $$W_0 = \C[p^{(n)}_\al,q^{(n)}_\al,u^{(n)}_i]_{\al \in
\De_+,i=1,\ldots,\el,n\geq 0}.$$ We identify it with $W_0$ defined in
the previous section by identifying
\begin{equation}    \label{coord}
p^{(n)}_\al \sim n!  a'_\al(-n-1), \quad \quad q^{(n)}_\al \sim n!
a^*_\al(-n), \quad \quad u^{(n)}_i \sim n! b_i(-n-1).
\end{equation}
The operator $\pa$ acts on $W_0$ as a derivation. There is a map $W_0
\otimes \C[t,t^{-1}] \arr \w$, which sends $P
\otimes t^n \in W_0[t,t^{-1}]$ to the corresponding local
functional, which we denote for simplicity by $\int P t^n$.

Using this map we can show that the space of local functionals
$\w_0$ is isomorphic to the quotient of $W_0$ by the subspace
of total derivatives and constants, cf. \cite{ff:im}. Thus, we can identify
the space of local functionals with $\w_0$, defined in the
previous section, as a linear space.

There is a Lie bracket on $\w_0$. It coincides with the well-known
Poisson structure, cf. \cite{ft,gd1}: $$\{ \int P,\int Q \} = - \sum_{1\leq
i,j\leq l} (\al_i,\al_j) \int \frac{\delta P}{\delta u_i} \, \pa \,
\frac{\delta Q}{\delta u_j} + \sum_{\al\in\De_+} \int \frac{\delta
P}{\delta p_\al} \, \frac{\delta Q}{\delta q_\al} - \sum_{\al\in\De_+} \int
\frac{\delta P}{\delta q_\al} \, \frac{\delta Q}{\delta p_\al},$$ where
$\delta/\delta f$ denotes variational derivative with respect to $f$.

This Lie bracket is uniquely defined by the Poisson brackets between
Fourier components of $p_\al(t), q_\al(t)$, and $u_i(t)$ are those given by
the formulas:
\begin{equation}    \label{class}
\{ p_\al(n),q_\al(m) \} = \delta_{\al,\beta} \delta_{n,-m}, \quad \quad \{
u_i(n),u_j(m) \} = n (\al_i,\al_j) \delta_{n,-m}.
\end{equation}

These formulas show that the Lie algebra structure on $\w_0$ for
generic $\nu$ is a quantum deformation of the Poisson Lie algebra structure
on $\w$ defined above: for any elements $A,B \in \w_0$, the
commutator is given by $[A,B] = \nu \{ A,B \} + \nu^2(\ldots)$.

Now we will give a hamiltonian interpretation of the classical limits of
the operators $G^\nu_i, i=1,\ldots,\el$.

Introduce formally variables $\phi_i, i=1,\ldots,\el$, such that $\pa \phi_i
= u_i$, and hence $\pa e^{\phi_i} = u_i e^{\phi_i}$ \cite{ff:im}. For any
$\la = \sum_{i=1}^\el \la_i \al_i$, put $W_\la = W_0 \otimes e^{\laa}$, where
$\laa = \sum_{i=1}^\el \la_i \phi_i$. We identify it with $W_\la$ defined in
the previous section by identifying $v_\la$ with $e^{\laa}$.

We have an action of $\pa$ on $W_\la$. We denote by $\w_\la$ the
quotient of $W_\la[t,t^{-1}]$ by the subspace of total
derivatives. This space can be interpreted as the space of functionals of
the form $\int P(p^{(n)}_\al(t),q^{(n)}_\al(t),u^{(n)}_i(t);t)
e^{\laa(t)} dt$. We will use simpler notation $\int P e^{\laa}, P \in W_\la$
for such a functional.

For any $\int P \in \w_0, \int Q e^{\laa} \in \w_\la$ define
their bracket $$\{ \int P,\int Q e^{\laa} \} = \sum_{1\leq i,j\leq l}
(\al_i,\al_j) \int \frac{\delta P}{\delta u_i} \, \left[ \la_i Q e^{\laa} -
\pa \left( \frac{\delta Q}{\delta u_j} e^{\laa} \right) \right] +$$
$$\sum_{\al\in\De_+} \int \frac{\delta P}{\delta p_\al} \, \frac{\delta
Q}{\delta q_\al} e\, ^{\laa(t)} - \sum_{\al\in\De_+}
\int \frac{\delta P}{\delta q_\al} \, \frac{\delta Q}{\delta p_\al} \,
e^{\laa(t)}.$$ This bracket is uniquely defined by formulas \eqref{class}
and
\begin{equation}    \label{classvertex}
\{ \int u_i t^n,\int e^{\laa} t^m \} = (\la,\al_i) \int e^{\laa} t^{n+m}.
\end{equation}

The map $\{ \cdot,\cdot \}: \w_0 \times \w_{-\al_i} \arr
\w_{-\al_i}$  satisfies the Jacobi identity for any $\int P_1, \int
P_2 \in \w_0, \int Q e^{\laa} \in \w_\la$ \cite{ff:im}. Hence it
defines on $\w_\la$ a structure of module over the Lie algebra
$\w_0$.

For $P \in W_0$ denote by $\xi(P)$ an operator on $W_\la$ given by
\begin{equation}    \label{xi}
\sum_{1\leq i,j\leq l} (\al_i,\al_j) \left[ \la_j \frac{\delta P}{\delta
u_j} + \sum_{n=1}^\infty \left( \pa^{n+1}
\frac{\delta P}{\delta u_i} \right) \frac{\pa}{\pa u_j^{(n)}} \right] +
\end{equation}
$$\sum_{\al\in\De_+} \sum_{n=1}^\infty \left( \pa^n \frac{\delta P}{\delta
p_\al} \right) \frac{\pa}{\pa q_\al^{(n)}} - \sum_{\al\in\De_+}
\sum_{n=1}^\infty \left( \pa^n \frac{\delta P}{\delta q_\al} \right)
\frac{\pa}{\pa p_\al}.$$

Clearly, $[\xi(P),\pa] = 0$ for any $P$. The map $\xi: W_0 \arr \on{End}
W_0$ is the classical limit of the map $P \arr \int Y(P,z) dz$ for generic
$\nu$, cf. \cite{ff:kdv}.

The map $\xi: W_0 \arr W_0$ can be presented in the form $\xi' \circ d$,
where $\xi'$ is a homomorphism $\Omega/\on{Im} \pa \arr \on{Vect}^\pa$ from
the quotient of the space of one-forms $\Omega$ by the action of $\pa$ to
the space of $\pa$--invariant vector field on $\on{Spec} W_0$. The map
$\xi'$ is a quasi-Poisson structure on $\on{Spec} W_0$ in the sense of
Gelfand-Dickey-Dorfman \cite{gd1,gd2} (cf. also \cite{ff:kdv}).

We can also define for any $P e^{\laa} \in W_\la$ a map $\xi(P e^{\laa}):
W_0 \arr W_\la$ by formula \eqref{xi}, where the first term should be
replaced by $$\sum_{1\leq i,j\leq l} (\al_i,\al_j) \left[ -\la_j P e^{\laa}
+ \sum_{n=1}^\infty \left(
\pa^{n+1} \frac{\delta P}{\delta u_i} e^{\laa} \right) \frac{\pa}{\pa
u_j^{(n)}} \right].$$ We have: $$\{ \int P,\int Q e^{\laa} \} = \int \left[
\xi(P) \cdot Q e^{\laa} \right] = - \int \left[ \xi(Q e^{\laa}) \cdot Q
\right],$$ cf. \cite{ff:im,ff:kdv}.

Now consider the following element of $W_{-\al_i}$: $$\GG_i = \sum_{\beta
\in \De_+} P^R_{\al_i,\beta}(q) p_\beta e^{-\phi_i}, \quad \quad
i=1,\ldots,\el,$$ where $P^R_{\al_i,\beta}(q)$ is a polynomial in $q_\gamma,
\gamma \in \De_+$, obtained from $P^R_{\al_i,\beta}$ in formula \eqref{LR}
by replacing $x_\gamma$ by $q_\gamma$. In the same way as in
\cite{ff:im,ff:kdv} we obtain the following result.

\begin{lem}    \label{ident}
$G_i^0 = \xi(\GG_i)$ and $\ol{G}_i^0 = \{ \cdot,\int \GG_i \}$.
\end{lem}

Thus, the space $J_0(\g)$ can be considered as the space of local
functionals in $p_\al(t), q_\al(t)$, and $u_i(t)$, which commute with $\int
\GG_i, i=1,\ldots,\el$, with respect to the Poisson bracket $\{ \cdot,\cdot
\}$. In other words, $J_0(\g)$ is the space of {\em local integrals of
motion} of the system of hamiltonian equations, defined by the hamiltonian
$$H = \sum_{i=1}^\el \int \GG_i.$$

This system reads: $$\pa_\tau p_\al(z,\tau) = \{ p_\al,H \}, \quad \quad
\pa_\tau q_\al(z,\tau) = \{ q_\al(z,\tau),H \}, \quad \quad \pa_\tau
u_i(z,\tau) = \{ u_i,H \}.$$ Here $p_\al(z,\tau), q_\al(z,\tau)$, and
$u_i(z,\tau)$ are considered as delta-like functionals, which are equal to
the value of the corresponding function at the point $t$, depending on the
time $\tau$.

Using formulas \eqref{class} and \eqref{classvertex}, we obtain the
following equations: $$\pa_\tau p_\al = \sum_{j=1}^\el
\sum_{\beta \in \De_+} \frac{\pa P^R_{\al_j,\beta}}{\pa q_\al} p_\beta
e^{-\phi_j}, \quad \quad \al \in \De_+,$$
\begin{equation}    \label{new}
\pa_\tau q_\al = - \sum_{j=1}^\el P^R_{\al_j,\al} e^{-\phi_j}, \quad \quad
\al \in \De_+,
\end{equation}
$$\pa_\tau \pa_z \phi_i = - \sum_{j=1}^\el (\al_j,\al_j) \sum_{\beta
\in \De_+} P^R_{\al_j,\beta} p_\beta e^{-\phi_j}, \quad \quad
i=1,\ldots,\el.$$

\smallskip
\noindent{\bf Example.} In the case of $\sw_2$ these equations read:
$$\pa_\tau p(z,\tau) = 0, \quad \quad \pa_\tau q(z,\tau) =
- e^{-\phi}, \quad \quad \pa_\tau \pa_z \phi(\tau,t) = p e^{-\phi}.$$
\qed

The system \eqref{new} should be compared to the system of Toda equations
associated to $\g$, which reads $$\pa_\tau u_i(z,\tau) = \sum_{j=1}^\el
(\al_i,\al_j) e^{-\phi_j}, \quad \quad i=1,\ldots,\el.$$ These
equations are non-local in $u_i(t)$, but possess local integrals of
motion. The corresponding algebra of integrals of motion is the classical
$\W$--algebra associated to $\g$ \cite{ff:im}.

\begin{rem}
The equations \eqref{new} imply that $$\pa_z \left( u_i -
\sum_{\al\in\De_+} (\al_i,\al) p_\al q_\al \right) = 0, \quad \quad
i=1,\ldots,\el.$$ Therefore we can put: $$u_i = U_i =
\sum_{\al\in\De_+} (\al_i,\al) p_\al q_\al$$ and eliminate $\phi_i,
i=1,\ldots,\el$ from the system \eqref{new} (recall that $u_i=\pa_z
\phi_i$). We then obtain the following system of equations on the
functions $p_\al, q_\al, \al \in \De_+$: $$\pa_\tau p_\al =
\sum_{j=1}^\el \sum_{\beta \in \De_+} p_\beta \frac{\pa
P^R_{\al_j,\beta}}{\pa q_\beta} e^{\int^z U_j dz}, \quad \quad \al \in
\De_+,$$ $$\pa_\tau q_\al = - \sum_{j=1}^\el P^R_{\al_j,\al} e^{\int^z
U_j dz}, \quad \quad \al \in \De_+.$$

Note that if we put: $p_{\al_i} = 1, i=1,\ldots,\el; p_\al = 0, \al
\neq \al_i, i=1,\ldots,\el$, in the resulting system, we will obtain a
system of equations on the functions $q_{\al_i}, i=1,\ldots,\el$,
which is equivalent to the Toda system. This is a version of
Drinfeld-Sokolov reduction.\qed
\end{rem}

We will now show that local integrals of motion of the system \eqref{new}
form the classical limit of the affine algebra $\G$.

Recall that the Lie group $N_+$ can be identified with the big cell of the
flag manifold $G/B_-$. Hence the Lie algebra $\g$ acts on $N_+$ from the
left by vector fields. We can write down formulas for these vector fields
in terms of the coordinates $x_\al, \al \in \De_+$: the action of elements
$e_\al$ is given by the vector fields $e_\al^L$, cf. formula \eqref{LR},
and $$h_i = - \sum_{\al\in\De_+} (\al_i,\al) x_\al
\frac{\pa}{\pa x_\al}, \quad \quad i=1,\ldots,\el,$$
\begin{equation}    \label{fi}
f_\al = \sum_{\beta \in \De_+} Q_{\al,\beta} \frac{\pa}{\pa x_\beta}, \quad
\quad \al \in \De_+,
\end{equation}
where $Q_{\al,\beta}$ is a certain polynomial of degree $\al+\beta$.

Introduce elements $E_i, F_i$, and $H_i, i=1,\ldots,\el$ of
$W_0$ as follows: $$H_i = u_i - \sum_{\al\in\De_+}
(\al_i,\al) p_\al q_\al,$$ $$E_i =
\sum_{\beta\in\De_+} P^L_{\al_i,\beta}(q) p_\beta, \quad \quad F_i =
\sum_{\beta\in\De_+} Q_{\al_i,\beta}(q) p_\beta + u_i q_{\al_i} +
\frac{2}{(\al_i,\al_i)} \pa q_{\al_i}.$$ Here
$P^L_{\al,\beta}(q)$ and $Q_{\al,\beta}(q)$ are obtained from the
polynomials $P_{\al,\beta}$ and $Q_{\al,\beta}$ given by \eqref{LR} and
\eqref{fi}, respectively, by replacing $x_\gamma, \gamma \in \De_+$, by
$q_\gamma$.

Now we define elements $E_\al$ and $F_\al$ for all other $\al \in \De_+$ by
induction. If $[e_\beta,e_\gamma] = e_{\beta+\gamma}$ in $\g$, then we put
$E_{\beta+\gamma} = \xi(E_\beta) \cdot E_\gamma$, where $\xi(\cdot)$ is
defined by formula \eqref{xi}, and analogously for $F_\al$.

\begin{thm}    \label{dp}
The space $K_0(\g)$ is the algebra of differential polynomials in
$E_\al, F_\al, \al\in\De_+$, and $H_i, i=1,\ldots,\el$, $$K_0(\g) =
\C[E_\al^{(n)},H_i^{(n)},F_\al^{(n)}]_{\al\in\De_+;i=1,\ldots,\el;n\geq
0}.$$ The space $J_0(\g)$ is the space of local functionals in $E_\al,
F_\al, \al\in\De_+$, and $H_i, i=1,\ldots,\el$, $$J_0(\g) =
K_0(\g)[t,t^{-1}] / (\on{Im} \pa \oplus \C).$$
\end{thm}

\begin{proof} One checks directly that the elements $E_\al,
F_\al$, and $H_j$ defined above lie in the kernel of the operators
$G_i^0, i=1,\ldots,\el$. Since $\pa$ commutes with $G_i^0$, the
derivatives of these elements also lie in the kernel. But $G_i^0 \cdot
(PQ) = (G_i^0 \cdot P) Q + P (G_i^0 \cdot Q)$. Hence the algebra
generated by these elements lies in $K_0(\g)$. But its character
equals the character of $K_0(\g)$, cf.  \thmref{zeroth}. Hence it
coincides with $K_0(\g)$.

The second part of the theorem follows from \thmref{zeroth}.
\end{proof}

Now consider the hyperplane $\G^*_1$ in the restricted dual space $\G^*$ to
the affine algebra $\G$, which consists of those linear functionals on $\G$
which are equal to $1$ on the central element $K \in \G$. This space has a
canonical Kirillov-Kostant Poisson structure.

If we choose a coordinate $t$ on the circle, we can identify $\G^*_1$ with
$\g[t,t^{-1}]$. The space ${\mc L}(\G)$ of local functionals on
$\g^*[t,t^{-1}]$ then becomes a Lie algebra with respect to the
Kirillov-Kostant bracket.

This bracket is uniquely defined by the brackets of the linear functionals
on $\G^*_1$. Such a functional is simply an element of $\g \otimes
\C[t,t^{-1}]$. The Kirillov-Kostant bracket of two such functionals, $A(n)
= A \otimes t^n$ and $B(m) = B \otimes t^m$ is given by $$\{ A(n),B(m) \} =
[A,B](n+m) + n \lef A,B \ri \delta_{n,-m}.$$

\begin{prop} $J_0(\g) \simeq {\mc L}(\G)$ as Lie algebras.
\end{prop}

\begin{proof}
We have to check that Fourier components of the differential polynomials
$E_\al, F_\al$, and $H_i$ have the same Poisson brackets as the
corresponding elements of ${\mc L}(\G)$. But this automatically follows
from our construction.
\end{proof}

Thus, we have shown that the space $J_0(\g)$ of local integrals of motion
of the system \eqref{new} is isomorphic to the space ${\mc L}(\G)$ of
local functionals on the dual space to the affine algebra $\G$. The space
${\mc L}(\G)$ is therefore the classical limit of the Lie algebra
$U_k(\G)_{\on{loc}}$ as $k \arr \infty$.

The space $K_0(\g)$ is the classical limit of the VOA $V_k$ as $k \arr
\infty$. We know from \thmref{zeroth} and \corref{gen} that
$K_\nu(\g)$ is also a VOA whose classical limit coincides with $K_0(\g)$.
Knowing the character of $K_\nu(\g)$ it is easy to show that $V_k$ is the
only possible quantum deformation of $K_0(\g)$. Finally, $k$ can be
computed from the commutator of the elements $H_i(n)$, which do not change
under deformation. Therefore $k=-h^\vee+\nu^{-1}$. This completes the
proof of \thmref{generic}.

This implies that $K_\nu(\g)$ coincides with $U_k(\G)_{\on{loc}}$ for
generic $\nu$. In particular, for any $k \neq -h^\vee$ there exists a free
field realization of the affine algebra $\G$ of level $k$, i.e. a
homomorphism $\G \arr U(\HH(\g) \oplus \h)_{\on{loc}}$, which sends $K$ to
$k$. It coincides with Wakimoto realization, which was first constructed
geometrically in \cite{ff:usp}.

\section{Integrals of motion of the deformed CFT}
We now want to add an extra operator $G_0^\nu$, corresponding to the
extra simple root $\al_0$ of $\G$, to the set $G_i^\nu,
i=1,\ldots,\el$. This operator will define certain deformations of
the conformal field theories associated to $\G$:
Wess-Zumino-Novikov-Witten model and generalized parafermions. We will
determine local integrals of motion of these deformations by an
analogue of the procedure from the previous section. The classical
limits of these integrals of motion will be shown to coincide with
local integrals of motion of the AKNS equation and its
generalizations.

\subsection{The case $\G=\su$.} In this case we have one operator
$G_1^\nu = \int a(z) V^\nu_{-\al}(z) dz$. It is natural to put
$$G_0^\nu = \int a^*(z) V^\nu_\al(z) dz.$$ The operator $G^\nu_0: W^\nu_0
\arr W^\nu_\al$ commutes with the action of $\pa$ and hence defines an
operator $\ol{G}_0^\nu: \w^\nu_0 \arr \w^\nu_\al$. We can now define a
Lie algebra $J_\nu(\su)$ as $$J_\nu(\su) = \on{Ker}_{\w^\nu_0}
\ol{G}_0^\nu \bigcap \on{Ker}_{\w^\nu_0} \ol{G}_1^\nu.$$

For $\la \in \C$, denote by $M^*_{\la,k}$ the module contragradient to the
Verma module over $\su$ with highest weight $\la$. There is a unique {\em
primary field} (or vertex operator) of weight $\la$ (i.e. spin $\la/2$)
$\Phi^\mu_\la(z): V_k \arr M^*_{\la,k}, \mu = \la,\la-2,\ldots$. Recall
that, by definition, the action of $\su$ on the components of this field by
commutator coincides with its action on the homogeneous components of the
evaluation representation corresponding to contragradient Verma module of
highest weight $\la$ over $\sw_2$, see, e.g., \cite{ff:lmp}.

We know that $\on{Ker}_{\w_0} G^\nu_1 = V_k$ with $k=-2+\nu^{-1}$ for
generic $\nu$, according to \thmref{generic}. Also $W^\nu_\al \simeq
M^*_{\al,k}$ as an $\su$--module for generic $\nu$. The restriction of the
fields $a^*(z)^m V^\nu_\la(z), m=0,1,\ldots: W^\nu_0 \arr W^\nu_\la$ to
$V_k$ is the primary field of weight $\la$, cf. \cite{ff:lmp}.  In
particular, $a^*(z) V^\nu_\al$ is the component $\Phi^0_2(z)$ of the
primary field corresponding to the adjoint representation of $\sw_2$ (we
can identify $\al$ with $2$).

Hence for generic $\nu$ elements of the space $J_\nu(\su)$ can be
interpreted as vectors $P \in V_k$, such that $\int \Phi^0_2(z) dz
\cdot P = \pa Q$ for some $Q \in W^\nu_\al$. For such a $P$, $\int
P(z) dz \in U_k(\G)\lo$ commutes with $\int \Phi^0_2(z) dz$. Following
Zamolodchikov \cite{zam}, we can interpret elements of $J_\nu(\su)$ as
local integrals of motion (in the first order of perturbation theory)
of the deformation of conformal field theory of $\su$ by the field
$\Phi^0_2(z)$. The space $J_\nu(\G)$ has a $\Z$--gradation, which is
obtained by subtracting $1$ from the $\Z$-gradation on $W_0^\nu$. We
will call an element of $J_\nu(\su)$ of degree $s$ an integral of
motion of spin $s$.

The field $\Phi^0_2(z)$ commutes with the homogeneous Heisenberg
subalgebra $\hh$ of $\su$. Hence it defines a primary field of the
coset model $\widehat{\sw}_2/\hh$ of arbitrary level $k \neq -2$,
which can be viewed as the analytic continuation of the parafermionic
theory \cite{fz,fateev} corresponding to the case of integral
$k$. This field is called the first thermal operator and is denoted by
$\ep_1(z)$. The corresponding deformation of the parafermionic theory
was studied by Fateev \cite{fateev}.

It is expected that for a positive integer $k$ the parafermionic
theory is equivalent to the $(k+1,k+2)$ minimal model of the conformal
field theory associated to the $\W$--algebra $\W(\sw_k)$. Hence we can
assume that this deformation is equivalent to the standard deformation
of the latter theory by the ``adjoint'' field \cite{EY,HM}, in which
spins of integrals of motion are presumably all positive integers,
which are not divisible by $k$. Thus we can expect the spins of
integrals of motion of the deformation of the parafermionic theory to
be the same, as proposed \cite{fateev}. This allowed Fateev to predict
the factorizable $S$--matrix of the deformed theory using the standard
bootstrap technique, see \cite{fateev}.

Since $\Phi^0_2(z)$ commutes with $\hh$, we have $U(\hh)\lo \subset
J_\nu(\su)$, and hence $J_\nu(\su)$ decomposes into the direct sum
$U(\hh)\lo \oplus J'_\nu(\su)$. By definition, elements of
$J'_\nu(\su)$ are local integrals of motion of the deformation of the
WZNW theory of level $k=-2+\nu^{-1}$ by $\int \Phi^0_2(z) dz$. Note
that the elements of $J'_\nu(\su)$ are also in one-to-one
correspondence with the integrals of motion of the deformation of the
parafermionic theory by $\int \ep_1(z) dz$.

\begin{conj}    \label{para}
For generic $\nu$ the space $J'_\nu(\su)$ is linearly spanned by
elements of all positive integral spins.
\end{conj}

We expect that when $k$ is an integer, the integral of motion of spins
divisible by $k$ indeed ``drops out'', in agreement with the
prediction of \cite{fateev}.

In the next section we will prove \conjref{para} in the classical
limit $\nu \arr 0$, and show that these integrals of motion are in
fact quantum deformations of the hamiltonians of the AKNS hierarchy.
We remark that the connection between the integrals of motion in the
deformations of the parafermionic theory and the non-linear
Schr\"odinger hierarchy (which is a reduction of the AKNS hierarchy)
has been previously dicussed by Schiff \cite{sch}.

Now we want to realize $J'_\nu(\su)$ as a cohomology group of a
complex $C_\nu^*(\su)$, which is constructed in the same way as the
complex $C_\nu^*(\g)$ in \secref{complex}.

We use the BGG resolution of $U_q(\su)$ with $q=\exp(\pi i \nu)$,
cf. \cite{ff:im}, Sects.~3.1, 4.4. Put $C_\nu^0(\su) = W^\nu_0$ and
$C_\nu^j(\su) = W^\nu_{2j} \oplus W^\nu_{-2j}, j>0$. The operators $a(m),
m\in\Z$, and $a^*(m), m \in \Z$, satisfy the Serre relations of
$\su$. Hence, according to \remref{ser}, the operators $G_1^\nu$ and
$G_0^\nu$ satisfy the $q$--Serre relations of $\su$ in the sense of
\lemref{qserre}. Thus we can define differentials of the complex
$C_\nu^*(\su)$ using the differentials of the BGG resolution of
$U_q(\su)$ in the same way as in \secref{complex}.

We have an analogue of \lemref{sing}: the differentials of the complex
$C_\nu^*(\su)$ are homogeneous operators, which commute with the action of
$\pa$. Therefore we can define the quotient complex $\cc^*(\su)$,
such that $\cc^j(\su) = C^*(\su)/\on{Im} \pa, j>0$, and
$\cc^0(\su) = C^0(\su)/(\on{Im} \pa \oplus \C)$. The $0$th
differential of $\cc^*(\su)$ is equal to $\ol{G}_1^\nu +
\ol{G}_0^\nu: \w_0^\nu \arr \w_{-\al}^\nu \oplus \W_\al^\nu$. Therefore
the space $J_0^\nu(\su)$ is the $0$th cohomology of the complex
$\cc^*(\su)$.

In the next section we will compute the cohomology of this complex in
the classical limit $\nu \arr 0$. But first we will define an analogue
of the space $J_\nu(\su)$ for an arbitrary affine algebra $\G$.

\subsection{General case.}    \label{imcft}
As an analogue of the operator $G_0^\nu$ for $\G$ we will take a
weight $0$ component from the vertex operator $\Phi_{\on{adj}}(z)$
corresponding to the adjoint representation of $\g$. The adjoint
representation is isomorphic to $\g$ itself, so its highest weight is
$$\al_{\on{max}} = -\al_0 = \sum_{i=1}^\el a_i \al_i,$$ where $\al_0$
is the weight corresponding to the $0$th simple root of $\G$. The
weight $0$ component is isomorphic to the Cartan subalgebra $\h$. More
generally, we can consider the vertex operator corresponding to the
contragradient Verma module $M_{\al_{\ma}}^*$ over $\g$ with highest
weight $\al_{\ma}$. This module contains the adjoint representation as
a submodule, and its component $M^*_{\al_{\ma}}(0)$ of weight $0$
contains $\h$.

To each vector $x \in M^*_{\al_{\ma}}(0)$ we can associate a field
$\Phi_{\on{adj}}^x(z): V_k \arr M^*_{-\al_0,k}$. For generic
$k=-h^\vee + \nu^{-1}$ this field has a bosonic realization $\Psi_x(z)
= P_x(a^*_\al(z)) V^\nu_{-\al_0}: W^\nu_0 \arr W^\nu_{-\al_0}$, where
$P_x(x_\al)$ is a polynomial in $x_\al, \al \in \De_+$, of weight
$-\al_{\on{max}}$ representing $x \in M^*_{\al_{\ma}}(0) \simeq
\C[x_\al]$.

Introduce the operator $$G_0^\nu = \int \Psi_x(z) dz:
W^\nu_0 \arr W^\nu_{-\al_0}.$$ It commutes with the action of $\pa$
and hence defines an operator $\ol{G}_0^\nu: \w^\nu_0 \arr
\w^\nu_{-\al_0}$. Define $$J_\nu(\G)_x = \bigcap_{i=0}^\el
\on{Ker}_{\w^\nu_0} G_i^\nu.$$ Note that $J_\nu(\G)_x$ is a Lie
subalgebra in $J^\nu(\g)$. It can be interpreted as the space of local
integrals of motion of the deformation of the WZNW model of level
$k=-h^\vee+\nu^{-1}$ by $\int \Phi_{\on{adj}}^x(z) dz$.

By definition, the field $\Phi_{\on{adj}}^x(z)$ commutes with the
action of the Heisenberg subalgebra $\hh$ of $\G$. Therefore we can
write $J_\nu(\G)_x = U(\hh)\lo \oplus J'_\nu(\G)_x$. The space
$J'_\nu(\G)_x$ can be considered as the space of local integrals of
motion of a deformation of the generalized parafermionic theory
associated to $\G$ (the $\g/\h$ coset model) by the field
$\Phi_{\on{adj}}^x(z)$.

Now we want to construct a complex $C_\nu^*(\G)_x$ whose cohomology
equals $J_\nu(\G)_x$, using the BGG resolution of $U_q(\G)$.  We put
$C_\nu^j(\G)_x = \oplus_{l(s)=j} W^\nu_{s(\rho)-\rho}$. In order to
define the differentials of the complex, we need to make sure that the
screening operators $G_i^\nu, i=0,\ldots,\el$, satisfy the $q$--Serre
relations of $\G$ in the sense of \lemref{qserre}. According to
\remref{ser}, this happens if the operators $e^R_i(n),
i=1,\ldots,\el$, and $P_x(m)$ satisfy the ordinary Serre
relations. Note that $P_x(m)$ automatically satisfies the relations
$[e^R_{i}(n_1),[e^R_{i}(n_2),P_x(m)]] = 0$ for all $i$. Therefore the
Serre relations hold for those $i$, for which the $i$th node of the
Dynkin diagram is not connected to the $0$th node. In addition, the
relations $[e^R_{i}(n),P_x(m)] = 0$ have to hold for all other $i$. In
other words, for those $i$'s the vector $x \in M^*_{\al_{\on{max}}}$
should be annihilated by the right action of $e_i$. In that case,
according to \remref{ser}, the operators $G_i^\nu, i=0,\ldots,\el$,
satisfy the $q$--Serre relations of $\G$ in the sense of
\lemref{qserre}.

Therefore we can construct differentials of the complex $C^*(\G)_x$
using the differentials of the BGG resolution of $\G$ in the same way
as in \secref{complex}. As before, the differentials of the complex
$C_\nu^*(\G)_x$ are homogeneous operators, which commute with the
action of $\pa$. Hence we can define the quotient complex
$\cc^*(\G)_x$, such that $\cc^j(\su) = C^j(\G)_x/ \on{Im} \pa, j>0$,
and $\cc^0(\su)_x = C^0(\G)_x/\on{Im} \pa \oplus \C$. The $0$th
differential of $\cc^*(\G)_x$ is equal to $\sum_{i=0}^\el
\ol{G}_i^\nu: \w_0^\nu \arr \oplus_{i=1}^\el
\w_{-\al_i}^\nu$. Therefore the space $J_0^\nu(\G)_x$ is the $0$th
cohomology of the complex $\cc^*(\su)_x$.

\begin{conj}    \label{paragen}
Suppose that $x \in M^*_{\al_{\on{max}}}(0)$ satisfies: $e_i^R \cdot x
= 0$ for all $i$, such that the $i$-th node of the Dynkin diagram is
not connected to the $0$th node. Then for generic $\nu$ the space
$J'_\nu(\G)_x$ has a linear basis, which consists of $\el$ elements of
each positive integral spin.
\end{conj}

In \secref{classgen} we will compute this cohomology of the complex
$C^*_\nu(\G)_x$ and prove \conjref{paragen} in the classical limit,
when $P_x = x_{\al_{\ma}}$ with respect to a special coordinate
system.

\section{Classical limit in the case of $\su$}

We consider now the new basis defined in \secref{class1}, which consists of
monomials in $a'(n)$, $a^*(n)$, and $b(n)$. This basis enables us to
identify the spaces $W^\nu_\la$ and $\w^\nu_\la$, respectively, for
different $\nu$, but consider the differentials of the complexes $C^*(\su)$
and $\cc^*(\su)$ as $\nu$--dependent, as in \secref{class1}.

In particular, we have for $G_1^\nu: W_0 \arr W_{-\al}$ and $G_0^\nu:
W_0 \arr W_\al$: $G_1^\nu = G_1^0 + \nu (\ldots), G_0^\nu = \nu G_0^0
+ \nu^2(\ldots)$. Put $G_1 = T_\al G_1^0: W_0 \arr W_0$, and $G_0 =
T_{-\al} G_0^0$.

\begin{lem}    \label{nilp2}
The operators $G_1$ and $G_0$ generate an action of the nilpotent
subalgebra $\N_+$ of $\su$ on $W_0$.
\end{lem}

This follows from the $q$--Serre relations to which the operators
$G_i^\nu$ satisfy in the limit $\nu \arr 0$. In \propref{isom} we will
give another proof of this fact. This implies

\begin{prop}    \label{cohiso}
The cohomology of the complex $C^*_0(\su)$ is isomorphic to the
cohomology of $\N_+$ with coefficients in $W_0$, $H^*(\N_+,W_0)$.
\end{prop}

\subsection{Separation of variables}    \label{sepa}
It is convenient to realize $W_0$ as
$\C[p^{(n)},q^{(n)},u^{(n)}]_{n\geq 0}$, where the variables $p^{(n)},
q^{(n)}$ and $u^{(n)}$ are defined by formula \eqref{coord}. In terms
of these variables we can write as in \lemref{ident}: $G_1^0 = \xi( -
p e^{-\phi}), G_0^0 = \xi( q e^\phi)$, and $\ol{G}_1^0 = \{ \cdot,-
\int p e^{-\phi} \}, \ol{G}_0^0 = \{ \cdot,\int q e^\phi \}$.

More explicitly, we have: $$G_1 = 2 \sum_{n=0}^\infty \pa^n(p
e^{-\phi}) \frac{\pa}{\pa u^{(n)}} + \sum_{n=0}^\infty \pa^n e^{-\phi}
\frac{\pa}{\pa q^{(n)}},$$ $$G_0 = 2 \sum_{n=0}^\infty \pa^n(q e^{\phi})
\frac{\pa}{\pa u^{(n)}} + \sum_{n=0}^\infty \pa^n e^{\phi} \frac{\pa}{\pa
p^{(n)}}.$$

By definition, $J_0(\su)$ is the kernel of the operator $\{ \cdot,\int
(-p e^{-\phi} + q e^{\phi}) \}$. Therefore we can interpret $J_0(\su)$
as the space of local integrals of motion of the system of equations:
$$\pa_\tau p(z,\tau) = e^\phi, \quad \quad \pa_\tau q(z,\tau) =
e^{-\phi},$$
$$\pa_\tau \pa_z \phi(z,\tau) = 2 q e^\phi + 2 p e^{-\phi}.$$

\begin{rem} These equations imply that $\pa_\tau (u - 2 pq) = 0.$
Hence we can put $u = 2pq$. Then we obtain the system
\begin{equation}    \label{new1}
\pa_\tau p(z,\tau) = e^{2\int^z pq dz}, \quad \quad \pa_\tau q(z,\tau)
= e^{-2\int^z pq dz}.
\end{equation}
If we identify $q(z,\tau) = \ol{p}(z,\tau)$ and replace $z$ by $iz$,
we obtain a non-local equation
\begin{equation}    \label{new2}
\pa_\tau p(z,\tau) = e^{- 2i \int^z |p|^2 dz}.
\end{equation}
The remarkable fact, which will be proved below, is that the integrals
of motion of the modified version of the AKNS hierarchy (resp.,
non-local Schr\"dinger hierarchy) are symmetries of equation
\eqref{new1} (resp., \eqref{new2}).\qed
\end{rem}

Recall the formula for the field $h(z)$ (see the Example at the end of
\secref{wreal}):
\begin{equation}    \label{hz}
h(z) = \sum_{n \in \Z} h_n z^{-n-1} = \frac{1}{\nu} b(z) - 2 :a(z)
a^*(z):.
\end{equation}
It is easy to check that all Fourier coefficients of this field (they
span the homogeneous Heisenberg subalgebra of $\wh{\sw}_2$) commute
with the screening operators $G_0^\nu$ and $G_1^\nu$ for all values of
$\nu$. Let us compute the classical limits of these Fourier
coefficients in terms of the variables $p^{(n)}, q^{(n)}, u^{(n)}$,
defined by formula \eqref{coord}. They will certainly commute with the
operators $G_0$ and $G_1$.

Set $v=u-2pq$. Then the limit of $\nu h_n, n<0$, as $\nu \arr 0$ is
the operator of multiplication by $\frac{1}{(-n-1)!} v^{(n)} =
\frac{1}{(-n-1)!} \pa^{-n-1} v$. On the other hand, a straightforward
calculation shows that $h_n, n>0$ equals $n! \pa_{n-1} +
\nu^2(\cdots)$, where
\begin{equation}    \label{pan}
\pa_n = 2 \frac{\pa}{\pa u^{(n)}} + 2
\sum_{m\geq 0} \left( \begin{array}{c} n+m+1 \\ m \end{array} \right)
\; \left[ p^{(m)} \frac{\pa}{\pa p^{(n+m+1)}} - q^{(m)}
\frac{\pa}{\pa q^{(n+m+1)}} \right], \quad \quad n \geq 0.
\end{equation}
Thus we obtain:
\begin{equation}    \label{comm}
[G_i^0,v^{(n)}]=0, \quad [G_i^0,\pa_n]=0, \quad \quad
i=0,1; n\geq 0,
\end{equation}
and
$$
[\pa_n,v^{(m)}] = 2 \delta_{n,m}.
$$
Finally, the leading term of $h_0$ as $\nu \arr 0$ equals
\begin{equation}    \label{h0}
h_0 = 2 \sum_{m\geq 0} \left[ p^{(m)} \frac{\pa}{\pa p^{(m)}} -
q^{(m)} \frac{\pa}{\pa q^{(m)}} \right].
\end{equation}

We find from the definition of $\pa_n$ the following relations:
$$
[\pa_n,\pa] = \pa_{n-1}, \quad n>0, \quad \quad [\pa_0,\pa] = h_0,
$$
This implies that the operator $\pa - \frac{1}{2} v h_0$ commutes with
all $\pa_n, n \geq 0$. Let us define the new variables
$\wt{p}^{(n)},\wt{q}^{(n)}, n \geq 0$, by the formulas
$$
\wt{p}^{(n)} = (\pa - \frac{1}{2} v h_0)^n p, \quad \quad \wt{q}^{(n)}
= (\pa - \frac{1}{2} v h_0)^n p.
$$
Then $\wt{p}^{(n)} = p^{(n)} + P_n, \wt{q}^{(n)} = q^{(n)} + Q_n$,
where $P_n, Q_n \in W_0$ lie in the ideal generated by $v^{(m)}, m\geq
0$, Since $\pa_m \cdot p = \pa_n \cdot q = 0$, we find that that
$\pa_m \wt{p}^{(n)} = \pa_m \wt{q}^{(n)} = 0$ for all $m\geq 0, n\geq
0$. Moreover, formulas \eqref{comm} imply that
$$
G_i \cdot (R(\wt{p}^{(n)},\wt{q}^{(n)}) P(v^{(n)})) = (G_i \cdot
R(\wt{p}^{(n)},\wt{q}^{(n)})) P(v^{(n)}).
$$

Therefore we can separate the variables, that is represent $W_0$ as
the tensor product $\C[\wt{p}^{(n)},\wt{q}^{(n)}]_{n\geq 0} \otimes
\C[v^{(n)}]_{n\geq 0}$, and $G_i$'s will act as derivations of the
factor $\wt{W}_0 = \C[\wt{p}^{(n)},\wt{q}^{(n)}]_{n\geq 0}$. It is
easy to find explicit expression for these derivations. Indeed, let us
identify $\C[\wt{p}^{(n)},\wt{q}^{(n)}]_{n\geq 0}$ with the quotient
of $W_0$ by the ideal generated by $v^{(n)}, n\geq 0$. Then $G_i \cdot
\wt{p}^{(n)}$ equals the projection of $G_i \cdot p^{(n)}n$ onto this
quotient, expressed as an element of
$\C[\wt{p}^{(n)},\wt{q}^{(n)}]_{n\geq 0}$. The same is true for $G_i
\cdot \wt{q}^{(n)}$. This way we obtain the following expressions for
$G_1$ and $G_0$ in the new variables:
$$G_1 = \sum_{n=0}^\infty B_n^- \frac{\pa}{\pa \wt{q}^{(n)}}, \quad
\quad G_0 = \sum_{n=0}^\infty B_n^+ \frac{\pa}{\pa \wt{p}^{(n)}},$$
where $B_n^\pm$ is defined recursively as follows: $B_0^\pm
= 1$, and
\begin{equation}    \label{rec}
B_n^\pm = \paa B_{n-1}^\pm \pm 2 \pq
B_{n-1}^\pm.
\end{equation}
Here we denote by $\paa$ the derivation of $\wt{W}_0$, such that $\paa
\wt{p}^{(n)} = \wt{p}^{(n+1)}, \paa \wt{q}^{(n)} = \wt{q}^{(n+1)}$.

\subsection{Isomorphism with $\C[N_+/H_+]$}
Let $N_+$ be the Lie group of $\N_+$. This is a prounipotent
proalgebraic Lie group, which is isomorphic to $\N_+$ via the
exponential map. Denote by $\hh_+$ the Lie subalgebra $\h \otimes
t\C[t]$ of $\N_+$ and by $H_+$ the corresponding subgroup of
$N_+$. The Lie algebra $\hh_+$ is called the homogeneous abelian Lie
subalgebra of $\N_+$. The Lie algebra $\N_+$ infinitesimally acts on
the homogeneous space $N_+/H_+$ from the left. On the other hand,
$N_+$ is isomorphic to the big cell of $B_- \backslash G$, and
therefore $N_+$ acts infinitesimally on $N_+$ from the right. The Lie
algebra $\hh_- = \h \otimes t^{-1} \C[t^{-1}]$ commutes with $\hh_+$,
and therefore it acts infinitesimally on $N_+/H_+$ from the
right. Thus, $\C[N_+/H_+]$ is an $\N_+$--module and an
$\hh_-$--module.

Let $E_1 = e \otimes 1$ and $E_0 = f \otimes t$ be the generators of
$\N_+$.

\begin{prop}    \label{isom}
There exists an isomorphism of rings
$\C[\wt{p}^{(n)},\wt{q}^{(n)}]_{n\geq 0} \simeq \C[N_+/H_+]$, under
which $E_i = - G_i, i=0,1$, and $\frac{1}{2} h_{-1} = \paa$.
\end{prop}

\begin{proof} We follow the same strategy as in the principal case
(see \cite{ff:kdv,Five}). Introduce the functions $\ol{p}$ and $\ol{q}$
on $N_+$ by the formulas
\begin{equation}
\ol{p}(K) = \frac{1}{2} (E_1,K h_{-1} K^{-1}), \quad \quad \ol{q}(K) =
- \frac{1}{2} (E_0,K h_{-1} K^{-1}), \quad \quad K \in N_+.
\end{equation}
These functions are invariant with respect to the right action of
$H_+$, and hence descend to $N_+/H_+$.

Next, we define the functions $$\ol{p}^{(n)} = (\frac{1}{2} h_{-1})^n
\cdot \ol{p}, \quad \quad \ol{q}^{(n)} = (\frac{1}{2} h_{-1})^n \cdot
\ol{q}, \quad \quad n>0,$$ on $N_+/H_+$. Define the homomorphism
$\C[\wt{p}^{(n)},\wt{q}^{(n)}]_{n\geq 0} \arr \C[N_+/H_+]$, which
sends $\ol{p}^{(n)}$ to $\wt{p}^{(n)}$, and $\ol{q}^{(n)}$ to
$\wt{q}^{(n)}$.

To prove that this homomorphism is injective, we have to show that the
functions $\ol{p}^{(n)}, \ol{q}^{(n)}, n\geq 0$ are algebraically
independent. We will do that by showing that the values of their
differentials at the identity coset $\bar{1} \in N_+/H_+$, are
linearly independent. Those are elements of the cotangent space to
$N_+/H_+$ at $\bar{1}$, which is canonically isomorphic to
$(\N_+/\hh_+)^*$. Using the invariant inner product on $\G$, we
identify $(\N_+/\hh_+)^*$ with $(\hh_+)^\perp \cap \N_- = \oplus_{j>0}
\n_+ \otimes t^{-j} \oplus \n_-) \otimes t^{-j+1}$.

Let us first show that the vectors $d \ol{p}|_{\bar{1}}$ and $d
\ol{q}|_{\bar{1}}$ generate $\n_+ \otimes t^{-1}$ and $\n_+ \otimes
t^{-1}$, respectively. For that it suffices to check that $E_1 \cdot
\ol{q} \neq 0, E_1 \cdot \ol{p} = 0, E_0 \cdot \ol{q} = 0, E_0 \cdot
\ol{p} \neq 0$. We find
\begin{align}
(E_1 \cdot \ol{q})(K) &= \frac{1}{2} (E_0,[E_1,K h_{-1} K^{-1}])
\notag = \frac{1}{2} ([E_0,E_1],K h_{-1} K^{-1}]) \\ & = - \frac{1}{2}
(h_1,K h_{-1} K^{-1}]) = - \frac{1}{2} (h_1,K h_{-1}
K^{-1}]) = - 1. \label{neg}
\end{align}
Likewise, we check that $E_1 \cdot \ol{p} = 0, E_0 \cdot \ol{q} = 0,
E_0 \cdot \ol{p} = -1$.

Now, by definition, $$d \ol{p}^{(m+1)}|_{\bar{1}} = \frac{1}{2}
\on{ad} h_{-1} \cdot d\ol{p}^{(m)}|_{\bar{1}}, \quad \quad d
\ol{q}^{(m+1)}|_{\bar{1}} = \frac{1}{2} \on{ad} h_{-1} \cdot
d\ol{q}^{(m)}|_{\bar{1}}.$$ Hence the vectors $d
\wt{u}^{(m)}_i|_{\bar{1}}$ are all linearly independent. But $\on{ad}
h_{-1}: \n_\pm \otimes t^{-j} \arr \n_\pm \otimes t^{-j-1}$ is an
isomorphism for all $j>0$. Therefore $d \ol{p}^{(n)}|_{\bar{1}}, d
\ol{q}^{(n)}|_{\bar{1}}, n\geq 0$ are linearly independent, and so the
functions $\ol{p}^{(n)}, \ol{q}^{(n)}, n\geq 0$, are algebraically
independent. Therefore our homomorphism
$\C[\wt{p}^{(n)},\wt{q}^{(n)}]_{n\geq 0} \arr \C[N_+/H_+]$ is
injective.

To prove that it is an isomorphism, we compute the characters of both
spaces with respect to the bigradation by the integers and the root
lattice. The result is
$$
\prod_{n=1}^\infty (1-q^n u^{-2})^{-1} (1-q^{n-1} u^2)^{-1}
$$
for both spaces, and so the above homomorphism is indeed an
isomorphism. Moreover, by construction, the operators $\frac{1}{2}
h_{-1}$ and $\paa$ get identified under this isomorphism.

It remains to show that the formula for $E_i$ in terms of the
coordinates $\ol{p}^{(n)}, \ol{q}^{(n)}$ coincides with the formula for $G_i,
i=0,1$. Let us show that for $E_1$. We already know that $E_1 \cdot
\ol{p} = -1, E_1 \cdot \ol{q} = 0$. On the other hand, in the same way
as in \cite{ff:kdv,Five} we obtain the relation
\begin{equation}    \label{recu}
[E_1,\frac{1}{2} h_{-1}] = - \frac{1}{2} f_{\al_1}(h_{-1}) \cdot E_1,
\end{equation}
where $f_{\al_1}(h_{-1})$ is the function on $N_+/H_+$, which equals
$(\al_1,K h_{-1} K^{-1})$ at $K \in N_+/H_+$. We claim that this
function equals $4 \pq$. To see this, note that $f_{\al_1}(h_{-1})$
has to be proportional to $\wt{p} \wt{q}$ by degree considerations. To
find the coefficient of proportionality, we apply to
$f_{\al_1}(h_{-1})$ the operator $E_1 E_0$. Our previous computations
show that $E_1 E_0 \cdot \pq = 1$, while we obtain in a similar
fashion: $E_1 E_0 \cdot f_{\al_1}(h_{-1}) = 4$.

Now formula \eqref{recu} gives us the relation
$$
[E_1,\paa] = - 2 \pq E_1.
$$
Writing
$$
E_1 = - \sum_{n\geq 0} \left( \wt{B}_n \frac{\pa}{\pa \wt{q}^{(n)}} +
\wt{B}'_n \frac{\pa}{\pa \wt{p}^{(n)}} \right),
$$
we find the recurrence relations on $\wt{B}_n, \wt{B}'_n$:
$$
\wt{B}_{n+1} = - 2 \wt{p} \; \wt{q} \wt{B}_n + \paa \wt{B}_n
$$
with the initial conditions $\wt{B}_0 = 1, \wt{B}'_0 = 0$. Comparing
this with formula \eqref{rec}, we obtain that $E_1=-G_1$. In the same
way we show that $E_0 = - G_0$. This completes the proof.
\end{proof}

\subsection{The cohomology}    \label{thecoh}
We can now compute the cohomology of the complex $C^*_0(\su)$.

\begin{prop}    \label{wedge}
The cohomology of the complex $C^*_0(\su)$ is isomorphic to
$\bigwedge^*(\hh_+^*) \otimes \C[v^{(n)}]_{n\geq 0}$.
\end{prop}

\begin{proof} According to the results of the previous subsection, as
an $\N_+$--module, $W_0 \simeq \C[N_+/H_+] \otimes \C[v^{(n)}]_{n\geq
0}$, where $\N_+$ acts trivially on the second factor. Therefore the
cohomology of the complex $C^*_0(\su)$, which by \propref{cohiso} is
isomorphic to $H^*(\N_+,W_0)$, equals $H^*(\N_+,N_+/H_+) \otimes
\C[v^{(n)}]_{n\geq 0}$. By Shapiro's lemma (see \cite{ff:kdv}),
$H^*(\N_+,N_+/H_+) \simeq \bigwedge^*(\hh_+^*)$.
\end{proof}

Now we want to compute the cohomology of the double complex
$\wh{C}^*_0(\su)$.

Denote
$$
\pa_v = \sum_{n\geq 0} v^{(n+1)} \frac{\pa}{\pa v^{(n)}}.
$$
The results of \secref{sepa} imply that 
$$
\pa = \frac{1}{2} h_{-1} + \pa_v + \frac{1}{2} v h_0.
$$
By \propref{wedge}, $h_n, n \leq 0$, act trivially on the cohomology
of the complex $C^*_0(\su)$. Hence $\pa$ acts on a representative of a
cohomology class $\omega \otimes P, \omega \in \bigwedge^*(\hh_+^*), P
\in \C[v^{(n)}]_{n\geq 0}$, as follows: $\pa \cdot \omega \otimes P =
\omega \otimes \pa_v P$. This leads to the following result.

\begin{prop}    \label{wedge0}
The $j$th cohomology of the complex $\wh{C}^*_0(\su)$ is isomorphic to
$\bigwedge^{j+1}(\hh_+^*) \oplus \bigwedge^j(\hh_+^*) \otimes
(\C[v^{(n)}]_{n\geq 0}/\on{Im} \pa_v)$.
\end{prop}

In particular, the zeroth cohomology is isomorphic to $\hh^* \oplus
\C[v^{(n)}]_{n\geq 0}/\on{Im} \pa_v$. The second summand consists of
classes of the form $\int P$, where $P \in \C[v^{(n)}]_{n\geq 0}$. The
classes corresponding to elements of the first summand are constructed
as follows.

Let $X \in W_{-\al} \oplus W_\al$ be a representative of a first
cohomology class of the form $\omega \otimes 1 \in \hh_+^* \otimes
\C[v^{(n)}]_{n\geq 0}$. Then $\pa X = 0$ in the cohomology, and hence
there exists $\wt{X} \in W_0$, such that $\delta^0 \wt{X} = X$ (we
recall that $\delta_0 = G_1^0 + G_0^0$). Here $\wt{X}$ is defined only
up to an element of $\C[v^{(n)}]_{n\geq 0}$, but there is a
representative of $h_0$--weight $0$, which lies in
$\C[\wt{p}^{(n)},\wt{q}^{(n)}]_{n\geq 0}$, and is unique up to a total
$\paa$--derivative. Note that $\wt{X}$ can not lie in the image of
$\pa$ (and hence $\paa$), because otherwise $X$ would also be in the
image of $\pa$. Therefore $\int h \neq 0$ defines a zeroth cohomology
class of $\wh{C}^*_0(\su)$.

Now let $\wt{X}_n, n<0$, be the element of $\wt{W}_0$, which
corresponds to the first cohomology class $h_{-n}^* \otimes 1$ of
$C^*_0(\su)$ via the isomorphism of \propref{wedge0}. Let $\eta_n =
\xi(\wt{X}_n)$. This is a derivation of $W_0$, which commutes with
$\pa$. It also has the following commutation relations with $G_1$ and
$G_0$:
$$[G_1,\eta_n] = - 2 \frac{\delta \wt{X}_n}{\delta u} G_1, \quad
[G_0,\eta_n] = 2 \frac{\delta \wt{X}_n}{\delta u} G_0.$$

Following \cite{ff:kdv} it is easy to describe all derivations $\xi$
of $W_0$, which commute with $\pa$ and have commutation relations of
the form
$$
[G_i,\xi] = (-1)^i f_\xi G_i,
$$
for some $f_\xi \in W_0$, with $G_i, i=0,1$.

\begin{lem} The vector space of such derivations is the direct sum of
$\hh_-$ and the space of derivations of the form $$\sum_{n\geq 0}
(\pa^{n+1}) P \frac{\pa}{\pa v^{(n)}} + \frac{1}{2} P h_0, \quad \quad P
\in \C[v^{(n)}]_{n\geq 0}.$$
\end{lem}

\begin{cor} $$\eta_n = \al_n h_{-n} + \sum_{n\geq 0} \pa^{n+1} P
\frac{\pa}{\pa v^{(n)}} + \frac{1}{2} P h_0$$ for some $\al_n \in
\C^\times, P \in \C[v^{(n)}]_{n\geq 0}$. In particular, $\{ \int
\wt{X}_n,\int \wt{X}_m \} = [\eta_n,\eta_m] = 0$ for all $n,m>0$.
\end{cor}

The Corollary means that the action of $\hh_-$ on $W_0$ is
hamiltonian.

\subsection{The zero curvature formalism}    \label{zcurv}

The evolutionary derivations of $\C[\wt{p}^{(n)},\wt{q}^{(n)}]$ coming
from the action of $\h_-$ can be written down explicitly in the zero
curvature form. This is explained in detail in \cite{EF1,Five} in the
case of the KdV and mKdV hierarchies, and the results carry over
directly to our case. The zero curvature equation corresponding to the
element $\frac{1}{2} h_{-n}, n>0$, of $\h_-$ reads:
\begin{equation}    \label{zc}
[\pa_z + (K \frac{1}{2}  h_{-1} K^{-1})_-,\pa_{\tau_n} + (K
\frac{1}{2}  h_{-n} K^{-1})_-] = 0, \quad \quad K \in N_+/H_+.
\end{equation}
Here for $A \in \g$, $A_-$ stands for the projection of $A$ onto the
${\mathfrak b}_-$ part of $\G = {\mathfrak b}_- \oplus
\N_+$. According to the computations made in the proof of
\propref{isom}, we have:
$$
(K \frac{1}{2} h_{-1} K^{-1})_- = \begin{pmatrix} pq + \frac{1}{2}
t^{-1} & -q t^{-1} \\ p & -pq - \frac{1}{2} t^{-1} \end{pmatrix}.
$$
(to simplify notation, we remove the tildes from $p$ and $q$). On the
other hand, for each $n>0$, $(K \frac{1}{2} h_{-n} K^{-1})_-$ is a
matrix, whose entries are polynomials in $t^{-1}$ with coefficients in
$\C[p^{(n)},q^{(n)}]$. Therefore formula \eqref{zc} defines a
derivation of $\C[p^{(n)},q^{(n)}]$. The first of the equations, with
$n=1$, tells us that $\tau_1 = z$, as expected. Straightforward
calculation gives:
$$
(K \frac{1}{2} h_{-1} K^{-1})_- =
$$
$$
= \begin{pmatrix} -2pq' + 2p'q - 2
p^2q^2 + pq t^{-1} + \frac{1}{2}
t^{-2} & q' -q t^{-2} \\ p' + p t^{-1} & 2pq' - 2p'q + 2
p^2q^2-pq t^{-1} - \frac{1}{2} t^{-2}
\end{pmatrix}
$$
(here $p'$ stands for $\pa p$). Substituting this into formula
\eqref{zc} we obtain the equation, corresponding to $n=2$:
\begin{equation}    \label{zc1}
\pa_{\tau_2} p = p'' - 2 p^3 q^2  -2 p^2 q', \quad \quad
\pa_{\tau_2} q = -q'' + 2 q^3 p^2 - q^2 p'.
\end{equation}

Now recall from \thmref{dp} that the kernel of the screening operator
$G_1$ in $W_0$ is $\C[E^{(n)},H^{(n)},F^{(n)}]_{n\geq 0}$, where $E =
p, H = v = u - 2pq, F = -pq^2 + uq + q'$. But $G_1^0$ commutes with
$v^{(n)} = H^{(n)}$ and preserves $\wt{W}_0 =
\C[\wt{p}^{(n)},\wt{q}^{(n)}]$. Therefore we obtain that the kernel of
$G_1^0$ in $\wt{W}_0$ equals $\C[\wt{E}^{(n)},\wt{F}^{(n)}]_{n\geq
0}$, where
$$
\wt{E} = \wt{p}, \quad \quad \wt{F} = \wt{p} \; \wt{q}^2 + \wt{q}'.
$$

The derivations $\pa_{\tau_n}$ commute with $G^0_1$ and hence define
evolutionary derivations of $\C[\wt{E}^{(n)},\wt{F}^{(n)}]_{n\geq 0}$,
which we denote by the same symbols. In particular, we find from
formula \eqref{zc1} the following formula for $\pa_{\tau_2}$ (we
again omit tildes to simplify notation):
\begin{equation*}    \label{zc2}
\pa_{\tau_2} E = E'' - 2 E^2 F, \quad \quad
\pa_{\tau_2} F = - F'' + 2 F^2 E.
\end{equation*}
This is the AKNS equation. Its reduction obtained by identifying $F$
with $\ol{E}$ and replacing $\tau_2$ with $i \tau_2$ is the non-linear
Schr\"odinger (nlS) equation
$$
i \pa_{\tau_2} = E'' - 2 E |E|^2.
$$
Therefore the derivations $\pa_{\tau_n}$, acting on
$\C[\wt{E}^{(n)},\wt{F}^{(n)}]_{n\geq 0}$, define the AKNS
hierarchy. Because of that, it is natural to call the equation
\eqref{zc1}, the modified AKNS (mAKNS) equation,
and the hierarchy of the derivations $\pa_{\tau_n}$ of
$\C[\wt{E}^{(n)},\wt{F}^{(n)}]_{n\geq 0}$, the mAKNS hierarchy.

The AKNS hierarchy also has a zero curvature representation
\eqref{zc}, where now $A_-$ stands for the projection of $A \in \G$
onto the $\g[t^{-1}]$ part of the decomposition $\G = \g[t^{-1}]
\oplus (\g
\otimes t\C[t])$, and $K \in \exp (\g \otimes t\C[t])/H_+$. In
particular, we then have the well-known AKNS $L$--operator
$$
\pa_z + (K \frac{1}{2} h_{-1} K^{-1})_- = \pa_z + \begin{pmatrix}
\frac{1}{2} t^{-1} & F \\ E & - \frac{1}{2} t^{-1}
\end{pmatrix}.
$$

Finally, the non-local equation \eqref{new1} can also be written in
the zero curvature form: $$ \left[ \pa_z + \begin{pmatrix} pq +
\frac{1}{2} t^{-1} & -q t^{-1} \\ p & -pq - \frac{1}{2} t^{-1}
\end{pmatrix},\pa_\tau + \begin{pmatrix} 0 & e^{-2\int^z pq dz} \\ t
e^{2\int^z pq dz} & 0 \end{pmatrix} \right] = 0.$$

The equations of the mAKNS hierarchy are symmetries of this equation.

\section{Classical limit for an arbitrary $\G$}    \label{classgen}

The results of this section can be generalized to the case of an
arbitrary (non-twisted) affine algebra $\g$. In this case $N_+/H_+$
also carries a good system of coordinates, in which the infinitesimal
action of $\hh_-$ becomes a set of evolutionary derivations. They
define a completely integrable system, which is an analogue of the
AKNS hierarchy. Moreover, with respect to these coordinates, the
generators $E_i, i=0,\ldots,\el$, of $\N_+$ become the classical
limits of the screening operators $G^\nu_i$.

\subsection{Screening operators in a special coordinate system}

Let us fix generators $e_\al, f_\al$ of the one-dimensional subspaces
$\n_{\al} \subset \n_+$ and $\n_{-\al} \subset \n_-$, respectively,
such that $(e_\al,f_\al) = 1$. In that case $[e_\al,f_\al] = \al$,
where we identify $\h$ and $\h^*$ using the inner product
on $\g$, such that $(\al_{\on{max}},\al_{\on{max}})=2$.

First we introduce a special coordinate system $\{ x_\al \}_{\al \in
\De_+}$ on the finite-dimensional unipotent group $\ol{N}_+$, with
respect to which the right action of $\n_+$ becomes particularly
simple. Let $\rh$ be the element of $\h$, such that $(\rh,\al_i) = 1,
\forall i=1,\ldots,\el$. Denote by $x_\al$ the regular function on
$\NN_+$ defined by the formula
$$
x_\al(\ol{K}) = - (f_\al,\ol{K}^{-1} \rh \ol{K}), \quad \quad \ol{K}
\in \NN_+.
$$
Then
\begin{align*}
(e_i^R \cdot x_{\al_i})(\ol{K}) & = - (f_{\al_i},[\ol{K}^{-1} \rh
\ol{K},e_i]) = - ([e_{\al_i},f_{\al_i}],\ol{K}^{-1} \rh \ol{K}) \\ & =
-(\al_i,\ol{K}^{-1} \rh \ol{K}) = -1.
\end{align*}

Let us write
$$
[e_\al,e_\beta] = - c_{\al,\beta} e_{\al+\beta}.
$$
Then
$$
[e_i,f_\al] = c_{\al_i,\al-\al_i} \; f_{\al-\al_i}, \quad \quad \al
\neq \al_i.
$$
Now we obtain in the same way as above:
$$
e_i^R \cdot x_\al = c_{\al_i,\al-\al_i} \; x_{\al-\al_i}, \quad \quad
\al \neq \al_i.
$$
Hence
$$
e_i^R = - \frac{\pa}{\pa x_{\al_i}} + \sum_{\al \in \Delta_+}
c_{\al_i,\al} \; x_\al \frac{\pa}{\pa x_{\al+\al_i}}.
$$
Moreover, we find that
\begin{align}    \label{eal}
e_\al^R &= - (\rh,\al) \frac{\pa}{\pa x_{\al}} + \sum_{\beta \in
\Delta_+} c_{\al,\beta} x_\beta \frac{\pa}{\pa x_{\al+\beta}}, \\
\notag h_i^L &= - \sum_{\al \in \De_+} (\al_i,\al) x_\al \frac{\pa}{\pa
x_\al}.
\end{align}

We use the above coordinates to construct the screening operators, as
in \secref{wreal}. Their classical limits as $\nu \arr 0$ are
\begin{equation}
\ol{G}^0_i = \left\{ \cdot,\int \left( -p_{\al_i} + \sum_{\al\in \De_+}
c_{\al_i,\al} \; q_\al p_{\al+\al_i} \right) e^{-\phi_i} \right\}
\end{equation}

Next we define the $0$th screening operator. According to
\secref{imcft} it has the form $P_x(a^*_\al(z)) V^\nu_{-\al_0}$, where
$P_x(x_\al) \in \C[\ol{N}_+]$ represents an element of
$M^*_{\al_{\ma}} \simeq \C[\ol{N}_+]$ of weight $0$, which satisfies
the conditions of \conjref{paragen}. It is straightforward to check
that the element $x_{\al_{\ma}}$ satisfies these conditions. In the
classical limit we obtain the following formula for the $0$th
screening operator:
$$
\ol{G}_0^0 = \left\{ \cdot,\int q_{\al_{\ma}} e^{-\phi_0} \right \}.
$$
The operators $\sum_{i=1}^\el \ol{G}^0_i$ and $\sum_{i=0}^\ell
\ol{G}^0_i$ define non-local Toda type equations. Here are the
explicit formulas for the second of these equations:
$$
\pa_\tau q_{\al_j} = e^{-\phi_j}, \quad \quad \pa_\tau q_\al =
\sum_{i=1}^\el c_{\al_i,\al-\al_i} q_{\al-\al_i} e^{-\phi_i}, \quad
\al\in\De_+\backslash\De_+^s,
$$
\begin{equation}    \label{nonloc1}
\pa_\tau p_\al = \sum_{i=1}^\el c_{\al_i,\al} e^{-\phi_i} +
e^{-\phi_0} \delta_{\al,\al_{\ma}}, \quad \quad \al\in\De_+,
\end{equation}
$$
\pa_\tau \pa_z \phi_j = \sum_{i=1}^\el (\al_i,\al_j) \left( p_{\al_i}
- \sum_{\al\in \De_+} c_{\al_i,\al} q_\al p_{\al+\al_i} \right)
e^{-\phi_i} - (\al_j,\al_0) q_{\al_{\ma}} e^{-\phi_0},
$$
$j=1,\ldots,\el.$

Let
\begin{equation}    \label{vi}
v_i = u_i - \sum_{\al\in\De_+} (\al_i,\al) q_\al p_\al.
\end{equation}
One finds that $G^0_i \cdot v_j = 0, \forall i,j$. Therefore we can
further reduce the system \eqref{nonloc1} by setting $u_i =
\sum_{\al\in\De_+} (\al_i,\al) q_\al p_\al$.

The classical screening operators give rise to derivations of the ring
of differential polynomials $W_0 =
\C[p_\al^{(n)},q_\al^{(n)},u_i^{(n)}]$, which we denote by $G_i,
i=0,\ldots,\ell$. We know that $G_i \cdot v_j = 0, \forall i,j,n$. In
the same way as in the case of $\su$ we can use this fact to separate
variables.

Each $h \in \h$ acts on $W_0$ in a natural way: $h \cdot
p_\al^{(n)} = \al(h) p_\al^{(n)}, h \cdot q_\al^{(n)} = - \al(h)
q_\al^{(n)}, h \cdot u_i^{(n)} = 0$. Denote by $\{ h^j
\}_{j=1,\ldots,\el}$, the dual basis to $\{ \al_i \}_{j=1,\ldots,\el}$
with respect to the normalized inner product on $\h$. Let $\paa = \pa
- \sum_{i=1}^\el v_i h^i$, and define new variables $\wt{p}_\al^{(n)}
= \paa^n p_\al, \wt{q}_\al^{(n)} = \paa^n q_\al$. In the same way as
in \secref{sepa}, one shows that $W_0 =
\C[\wt{p}_\al^{(n)},\wt{q}_\al^{(n)}] \otimes \C[v_i^{(n)}]$, and that
the derivations $G_i$ act along the first factor of this tensor
product.

Furthermore, we find, in the same way as in \secref{sepa}, the
following explicit formulas for the action of $G_i, i=1,\ldots,\el$,
on $\C[\wt{p}_\al^{(n)},\wt{q}_\al^{(n)}]$:
\begin{equation}    \label{Gi}
G_i = \sum_{n=0}^\infty \sum_{\al\in\De_+} \left( B_{\al,n}^{(i)}
\frac{\pa}{\pa \wt{q}^{(n)}_\al} + B_{-\al,n}^{(i)} \frac{\pa}{\pa
\wt{p}^{(n)}_\al} \right),
\end{equation}
where the polynomials $B_{\pm\al,n}$ are defined recursively as
follows:
\begin{align*}
B_{0,\al_j}^{(i)} &= \delta_{i,j}; \quad \quad B_{0,\al}^{(i)} = -
c_{\al_i,\al-\al_i} \wt{q}_{\al-\al_i}, \quad \al \in
\De_+\backslash\De_+^s, \\ B_{-\al,0}^{(i)} & = c_{\al_i,\al}
\wt{p}_{\al+\al_i}, \\ B_{\pm\al,n}^{(i)} &= \paa B_{\pm\al,n-1}^{(i)}
- U_i B_{n-1}^\pm, \quad \quad n>0,
\end{align*}
where
\begin{equation}    \label{Ui}
U_i = \sum_{\al\in\De_+} (\al_i,\al) \wt{p}_\al \wt{q}_\al.
\end{equation}
We also find a formula for $G_0$:
\begin{equation}    \label{G0}
G_0 = \sum_{n=0}^\infty B_n^{(0)} \frac{\pa}{\pa
\wt{p}^{(n)}_{\al_{\ma}}},
\end{equation}
where
$$
B_0^{(0)} = 1, \quad \quad B_n^{(0)} = \paa B_{n-1}^{(0)} - U_0
B_{n-1}^{(0)},
$$
with $U_0$ given by formula \eqref{Ui}.

\subsection{Isomorphism with $\C[N_+/H_+]$}

Now we can identify the ring of differential polynomials
$\C[\wt{p}_\al^{(n)},\wt{q}_\al^{(n)}]$ with $\C[N_+/H_+]$, where
$H_+$ is the subgroup of $N_+$, which is the image of the Lie algebra
$\h \otimes t \C[[t]]$ under the exponential map. Observe that the Lie
algebra $\h \otimes t^{-1} \C[t^{-1}]$ acts on $\C[N_+/H_+]$ from the
right. Denote by $\rhh$ the element $\rh \otimes t^{-1}$ of $\h
\otimes t^{-1} \C[t^{-1}]$. Let $E_i = e_i \otimes 1, i=1,\ldots,\el$,
and $E_0 = f_{\al_{\ma}} \otimes t$ be the generators of $\N_+$.

\begin{prop}    \label{isom1}
There exists an isomorphism of rings
$\C[\wt{p}^{(n)},\wt{q}^{(n)}]_{n\geq 0} \simeq \C[N_+/H_+]$, under
which $E_i = - G_i, i=0,\ldots,\el$, and $\rhh = \paa$.
\end{prop}

\begin{proof} The proof proceeds along the lines of the proof of
\propref{isom}. We introduce a system of coordinates on $N_+/H_+$ and
then find formulas for the action of $E_i = - G_i, i=0,\ldots,\el$,
and $\rhh$ in these coordinates.

Introduce the following regular functions $\ol{p}_\al$, $\ol{q}_\al$
on $N_+/H_+$:
\begin{align*}
\ol{p}_\al(K) &= - (e_\al \otimes 1,K \rhh K^{-1}), \quad \quad K \in
N_+/H_+, \\
\ol{q}_\al(K) &= - (f_\al \otimes t,K \rhh K^{-1}), \quad \quad K \in
N_+/H_+.
\end{align*}
Straightforward computation analogous to that made in the proof of
\propref{isom} gives for $i=1,\ldots,\el$:
\begin{align*}
E_i \cdot \ol{q}_\al &= - \delta_{i,j}; \quad \quad E_i \cdot
\ol{q}_\al = c_{\al_i,\al-\al_i} \; \ol{q}_{\al-\al_i}, \quad \al \in
\De_+\backslash\De_+^s, \\ E_i \cdot \ol{p}_\al & = - c_{\al_i,\al} \;
\ol{p}_{\al+\al_i},
\end{align*}
and
\begin{equation}    \label{E0}
E_0 \cdot \ol{q}_\al = 0, \quad \quad E_0 \cdot \ol{p}_\al =
(\al_{\ma},\rh) \delta_{\al,\al_{\ma}} -
c_{\al,\al_{\ma}-\al} q_{\al_{\ma}-\al}.
\end{equation}

Now let $\ol{p}^{(n)}_\al = (\rhh)^n \cdot \ol{p}^{(n)}_\al,
\ol{q}^{(n)}_\al = (\rhh)^n \cdot \ol{q}^{(n)}_\al, n\geq 0$. We show
in the same way as in the proof of \propref{isom} that these functions
are algebraically independent. Here we rely only on the fact that that
$\rh$ is a regular semi-simple element of $\g$.

Next we define a homomorphism from
$\C[\wt{p}_\al^{(n)},\wt{q}_\al^{(n)}]$ to $\C[N_+/H_+]$. It
intertwines the actions of $\paa$ and $\rhh$ and maps each
$\wt{q}_\al$ to $\ol{q}_\al$. It also maps each $\wt{p}_\al$ to a
polynomial in $\ol{p}_\al, \ol{q}_\al$ (also denoted by $\wt{p}_\al$),
such that
\begin{equation}    \label{change}
\ol{p}_\al = (\al,\rh) \wt{p}_\al - \sum_{\beta \in \Delta_+}
c_{\al,\beta} \ol{q}_\al \wt{p}_{\al+\beta}
\end{equation}
(cf. formula \eqref{eal}). It is clear that such polynomials exist,
are unique and that the above homomorphism is injective. To prove that
it is an isomorphism, we use the equality of characters of the two
spaces, as in the proof of \propref{isom}.

Finally, we need to show that the action of $E_i$ on $\C[N_+/A_+]$
coincides with the action of $-G_i$ on
$\C[\wt{p}_\al^{(n)},\wt{q}_\al^{(n)}]$. But the actions of $E_i$ and
$-G_i$ on $\ol{q}_\al = \wt{q}_\al$ coincide. The action of $E_i,
i=1,\ldots,\el$, on $\wt{p}_\al$ is given by the same formula as the
action on $\ol{p}_\al$. Comparing with formula \eqref{Gi} we again
find agreement between the actions of $E_i$ and $-G_i$. The same is
true for $E_0$ and $-G_0$ as formulas \eqref{E0}, \eqref{G0} and
\eqref{change} show.

We also have:
\begin{equation}    \label{commut}
[E_i,\rhh] = - f_{\al_i} \cdot E_i,
\end{equation}
where $f_{\al_i}(K) = (\al_i,K \rhh K^{-1}), K \in N_+/H_+$. It is
easy to see that
$$
f_{\al_i} = \sum_{\al\in\De_+} (\al_i,\al) p_\al q_\al.
$$
Formula \eqref{commut} gives us a recurrence relation on the
coefficients of the derivations $E_i$, which coincides with that on
the coefficients of the derivations $G_i$ given by formulas
\eqref{Gi}, \eqref{G0}. We have shown above that their first
coefficients differ by sign, and therefore we obtain that $E_i = -
G_i$. This completes the proof.
\end{proof}

\begin{rem} In the above isomorphism the element $\rh$ of $\h$ can be
replaced by any regular element of $\h$. In that case we need to
change accordingly the formulas for the coordinates $p_\al,
q_\al$. Moreover, with appropriate changes all results of this section
will remain true if we replace $\rh$ by any regular element of
$\h$.\qed
\end{rem}

\subsection{Integrable hierarchies}
Recall that the Lie algebra $\hh_- = \h \otimes t^{-1} \C[t^{-1}]$
acts on $\C[N_+/H_+]$ by derivations. Hence we obtain an infinite
hierarchy of commuting evolutionary (i.e., commuting with $\paa =
\rhh$) derivations on $\C[\wt{p}_\al^{(n)},\wt{q}_\al^{(n)}]$. We call
it the modified AKNS hierarchy associated to $\g$, or $\g$--{\em mAKNS
hierarchy} for shorthand.

The above derivations preserve the subring $\C[\ol{N}_+\backslash
N_+/H_+]$ of $\ol{N}_+$--invariants of $\C[N_+/H_+]$. This subring
equals the intersection of kernels of the operators $G_i,
i=1,\ldots,\el$, in $\C[\wt{p}_\al^{(n)},\wt{q}_\al^{(n)}]$. The
latter is isomorphic to the quotient of the ring $K_0(\g) =
\C[E_\al^{(n)},H_i^{(n)},F_\al^{(n)}]$ introduced in \thmref{dp} by
the ideal generated by $H_i^{(n)}$. Hence it is also isomorphic to a
ring of differential polynomials
$\C[\wt{E}_\al^{(n)},\wt{F}_\al^{(n)}]$, where $\wt{E}_\al^{(n)}$ and
$\wt{F}_\al^{(n)}$ are certain polynomials in
$\wt{p}_\beta^{(m)},\wt{q}_\beta^{(m)}$. The Lie algebra $\hh_-$ acts
on $\C[\wt{E}_\al^{(n)},\wt{F}_\al^{(n)}]$ by evolutionary
derivations. They form a hierarchy, which we call the $\g$--{\em AKNS
hierarchy}. It has been previously studied in the literature under the
name AKNS--D hierarchy (see \cite{Di}).

The equations of both hierarchies can be written in the zero curvature
form as follows. The equation of the $\g$--mAKNS hierarchy
corresponding to an element $y \in \hh_-$ reads
$$
[\pa_z + (K \rhh K^{-1})^-,\pa_{\tau} + (K y K^{-1})^-] = 0, \quad
\quad K \in N_+/H_+.
$$
where $A^-$ denotes the projection of $A \in \G$ onto the ${\mathfrak
b}_-$ part of the decomposition $\G = {\mathfrak b}_- \oplus \N_+$.

The equation of the $\g$--AKNS hierarchy corresponding to
$y \in \hh_-$ reads
$$
[\pa_z + (K \rhh K^{-1})_-,\pa_{\tau} + (K y K^{-1})_-] = 0.
$$
where $A_-$ stands for the projection of $A \in \G$ onto the
$\g[t^{-1}]$ part of the decomposition $\G = \g[t^{-1}] \oplus (\g
\otimes t\C[t])$, and $K \in \exp (\g \otimes t\C[t])/H_+$.
In the next subsection we will show that both hierarchies are
hamiltonian.

Note also that the non-local equations introduced above can also be
written in the zero-curvature form (in the case of $\su$ see the end
of \secref{zcurv}). The equations of the $\g$--mAKNS hierarchy are
symmetries of the non-local equation \eqref{nonloc1}.

\subsection{Cohomology computation}

In the same way as in the case of $\su$ we obtain the following
result.

\begin{prop}    \label{wedge1}
The cohomology of the complex $C^*_0(\G)$ is isomorphic to
$H^*(\N_+,W_0)$ and hence to $\bigwedge^*(\hh_+^*) \otimes
\C[v_i^{(n)}]_{i=1,\ldots,\el;n\geq 0}$.
\end{prop}

Now we can to compute the cohomology of the double complex
$\wh{C}^*_0(\su)$.

We have
$$
\pa = \rhh + \sum_{i=1}^\el v_i h^i + \pa_v,
$$
where
$$
\pa_v = \sum_{i=1}^\el \sum_{n\geq 0} v_i^{(n+1)} \frac{\pa}{\pa
v_i^{(n)}}.
$$
In the same way as in the case of $\su$ we obtain:

\begin{prop}    \label{wedge2}
The $j$th cohomology of the complex $\wh{C}^*_0(\G)$ is isomorphic to
$\bigwedge^{j+1}(\hh_+^*) \oplus \bigwedge^j(\hh_+^*) \otimes
(\C[v_i^{(n)}]/\on{Im} \pa_v)$.
\end{prop}

In particular, the zeroth cohomology is isomorphic to $\hh^* \oplus
\C[v_i^{(n)}]/\on{Im} \pa_v$. The classes corresponding to elements of
$\hh^*$ are constructed as in the case of $\su$.

Let $\{ h^i \}_{i=1,\ldots\el;n<0}$ be a basis of $\h$. Then $\{ h^i_n
\}_{i=1,\ldots,\el;n>0}$, where $h^i_n = h^i \otimes t^{-n}$ is a
basis of $\h \otimes t^{-1} \C[t^{-1}]$. Denote by $\wt{X}^i_n \in
W_0$, the representative of the zeroth cohomology class of
$\wh{C}^*_0(\su)$, which corresponds to the first cohomology class
$(h^i_{-n})^* \otimes 1$ of $C^*_0(\G)$ via the isomorphism of
\propref{wedge2}. Let $\eta^i_n = \xi(\wt{X}^i_n)$. This is a
derivation of $W_0$, which commutes with $\pa$ and has the commutation
relations with the operators $G_i$ of the form
\begin{equation}    \label{allsuch}
[G_i,\xi] = f^\xi_{\al_i} G_i, \quad \quad i=0,\ldots,\el,
\end{equation}
for some $f^\xi_{\al_i} \in W_0$ (in the case of $\su$ these relations
are given in \secref{thecoh}).

Following \cite{ff:kdv} it is easy to describe all such derivations
$\xi$ of $W_0$.

\begin{lem} The vector space of derivations $\xi$ of $W_0$ satisfying
commutation relations \eqref{allsuch} is the direct sum of $\hh_-$ and
the space of derivations of the form $$\xi(\{P_j\}) = \sum_{i=1}^\el
\sum_{n\geq 0} (\pa^{n+1} P_i) \frac{\pa}{\pa v_i^{(n)}} +
\sum_{i=1}^\el P_i h^i, \quad \quad P_j \in \C[v_i^{(n)}].$$
\end{lem}

\begin{cor} $$\eta^i_n = \al^i_n h^i_{-n} + \xi(\{P_j\})$$ for some
$\al^i_n \in \C^\times, P_j \in \C[v_i^{(n)}]$. Therefore $\{ \int
\wt{X}^i_n,\int \wt{X}^j_m \} = [\eta^i_n,\eta^j_m] = 0$ for all
$i,j,n,m$.
\end{cor}

The Corollary implies that the $\g$--mAKNS and $\g$--AKNS hierarchies
introduced above are completely integrable hamiltonian systems.

\end{document}